\documentclass{ba-kecs}
\title{Connectivity Structure of Systems}
\runningtitle{Connectivity structure of systems}
\author{Remco Bras\\Bachelor Thesis\\Department of Knowledge Engineering\\Maastricht University}
\runningauthor{Remco Bras}
\usepackage{amsthm,amsmath,amssymb}
\usepackage{graphicx}
\usepackage{xspace}
\newtheorem{lemma}{Lemma}
\newtheorem{theorem}{Theorem}
\newtheorem{observation}{Observation}
\theoremstyle{definition}
\newtheorem{definition}{Definition}
\theoremstyle{remark}
\newtheorem{example}{Example}
\newtheorem{corollary}{Corollary}
\newtheorem{remark}{Remark}
\newcommand{\charmat}[1]{\ensuremath{\lambda I - #1}}
\newcommand{\dcharmat}[1]{\ensuremath{\left|\charmat{#1}\right|}}
\newcommand{\stdnamedsys}[0]{\namedsys{S}{A}{B}{C}{D}}
\DeclareMathOperator{\adj}{adj}
\newcommand{\cadj}[1]{\ensuremath{\adj\left(\charmat{#1}\right)}}
\newcommand{\namedsys}[5]{\ensuremath{#1 = \left(#2,#3,#4,#5\right)}\xspace}
\newcommand{\transys}[2]{\ensuremath{\mathbb{T}\left(#1,#2\right)}}

\newcommand{\structuredspace}[2]{\ensuremath{\mathbb{S}_{\mathbb{R}}^{#1
      \times #2}}\xspace}
\newcommand{\namedsysstruct}[5]{\ensuremath{#1 =
    \left(#2,#3,#4,#5\right)^S}\xspace}
\newcommand{\stdnamedsysstruct}[0]{\namedsysstruct{\mathbf{S}}{A}{B}{C}{D}}
\newcommand{\dual}[1]{\ensuremath{\mathbb{D}\left(#1\right)}\xspace}
\newcommand{\dualstruct}[1]{\ensuremath{\mathbb{D}^S\left(#1\right)}\xspace}
\begin{document}
\maketitle
\begin{abstract}
In this paper, we consider to what degree the structure of a linear system is
determined by the system's input/output behavior.
The structure of a linear system is a directed graph where the
vertices represent the variables in the system and an edge \((x,y)\)
exists if \(x\) directly influences \(y\).
In a number of studies, researchers have attempted to identify such
structures  using input/output data.
Thus, our main aim is to consider to what degree the results of such
studies are valid.
We begin by showing that in many cases, applying a linear
transformation to a system will change the system's graph.
Furthermore, we show that even the graph's components and their
interactions are not determined by input/output behavior.
From these results, we conclude that without further assumptions, very
few aspects, if any, of a system's structure are determined by its
input/output relation.
We consider a number of such assumptions.
First, we show that for a number of parameterizations, we can
characterize when two systems have the same structure.
Second, in many applications, we can use domain knowledge to exclude certain
interactions.
In these cases, we can assume that a certain variable \(x\)
does not influence another variable \(y\).
We show that these assumptions cannot be sufficient to identify
a system's parameters using input/output data.
We conclude that identifying a system's structure
from input/output data may not be possible given only assumptions of
the form \(x\) does not influence \(y\).

\end{abstract}
\section{Introduction}
The aim of this paper is to consider to what degree the structure of a
linear system can be determined using input/output data.
By the structure of a linear system, we mean a graph in which the
vertices represent the variables in the system and an edge from one
variable to another exists if the second variable is directly
influenced by the first one. 

\subsection{Motivation}
This problem arose in a number of recent studies 
in which researchers have attempted to find the structure of a dynamical system using
input/output data.
For instance, researchers have attempted to find
structure in the brain using fMRI data.
In particular, Friston et al.\ \cite{friston} propose the method of dynamic causal modeling.
In this method, the activity of regions in the brain is modeled by a
bilinear system.
Each state variable \(z_i\) of this system represents the activity in a
particular region of the brain.
The change over time of this activity is given by
\eqref{eq:bilinear-state}.
\begin{equation}
\label{eq:bilinear-state}
  \dot{z} = (A+\sum_{i=1}^{n_u}u_iB_i)z+Cu
\end{equation}

Here, \(u\) is the vector of inputs \(u_i\) to the model and \(n_u\)
the number of such inputs.
In dynamic causal modeling, these inputs represent experimental
conditions.
For some of these inputs, the corresponding column in the matrix \(C\)
will be non-zero, implying that they affect state variables directly.
For other inputs, the corresponding matrix \(B_i\) is non-zero,
allowing these inputs to indirectly influence the state variables by
changing how these variables affect each other.

In addition to the bilinear model describing brain activity, Friston
et al.\ use what they call a ``forward model'' that describes how this
activity is measured.
This model depends on the particular measurement method used, such as
fMRI or EEG.
A forward model for fMRI measurements is given by Friston et
al.\ \cite{friston}.

To find structure in the brain, Friston et al.\ identify the parameters
of both the bilinear model and the forward model.
To do so, they assume the bilinear model has a particular form and the
parameters that appear in this form and the forward model have certain given prior probability densities.
Friston et al.\ then use the data to find the posterior densities of
the parameters.
These densities can be used to make inferences about the parameters.
For instance, by testing whether a particular entry of \(A\) is
non-zero, we can test if the data supports the hypothesis that one
variable influences another in the absence of input.

Goebel et al.\ \cite{goebel} propose a different approach based on vector
autoregressive (AR) models.
In these models, a vector time-series \(x_n\) is computed using its
own past values, as shown in \eqref{eq:vector-ar}.
Here, the integer \(p\) is called the order of the vector AR model.
The input \(u_n\) is a stochastic white noise input with a given
cross-covariance matrix.
Unlike Friston et al., Goebel et al.\ do not use inputs based on
experimental conditions.

\begin{equation}
  \label{eq:vector-ar}
  x_n = -\sum_{i=1}^pA_ix_{n-i} + u_n
\end{equation}

To quantify the influence of one variable on another, Goebel et
al.\ use the concept of Granger causality.
Goebel et al.\ distinguish two forms of this concept.
The first of these is directed influence from one time series \(x\) to a
time series \(y\).
We say that \(x\) causes \(y\) if we can better predict \(y_n\)
using past values of both \(x\) and \(y\), that is, the set
\(S_{x,y} = \{x_{n-1},y_{n-1},x_{n-2},y_{n-2},\cdots\}\), than using past values
of \(y\) alone, i.e. the set \(S_{y} = \{y_{n-1},y_{n-2},\cdots\}\).
The second form of Granger causality is instantaneous causality
between \(x\) and \(y\).
This form of causality occurs if we can better predict \(y\) from
\(S_{x,y} \cup \{x_n\}\) than from \(S_{x,y}\).
Goebel et al.\ note that though the first form of causality is
directed, the second is not.

To apply the concept of Granger causality to vector AR models, Goebel
et al.\ create three such models.
One model predicts \(x\) in terms of its past values and another model
performs a similar task for \(y\).
The third model uses past values of both \(x\) and \(y\) to predict
both \(x\) and \(y\).
Goebel et al.\ then use a number of measures based on the covariance
matrices of the noise vectors that appear in these models to quantify
the presence of Granger causality between \(x\) and \(y\).

The two papers we have considered so far focused on identifying
structure in the brain from fMRI data. 
A similar method can also be
applied in other fields.
An example of this is the application described in Hollanders' PhD
thesis \cite{hollanders}.
In this application, Hollanders considers time series of the
expression levels of the genes of a unicellular fungus.
Hollanders' goal was to identify how these expression levels influence
each other.
To do so, Hollanders considers two model classes.
One of these is the class of linear systems, the other a
generalization of linear systems called piecewise linear systems.
The state equation for one of Hollanders' linear systems is given by
\eqref{eq:continuous-lti}.
In this equation, \(x(t)\) is the vector of gene expressions at time
\(t\), \(u(t)\) is a vector of inputs and \(\xi_t\) is a vector of Gaussian
white noise.
The matrices \(A\) and \(B\) are constant.

\begin{equation}
  \label{eq:continuous-lti}
  \dot{x}(t) = Ax(t) + Bu(t) + \xi_t
\end{equation}

Hollanders' second type of system, a piecewise linear system, is given
by \eqref{eq:continuous-piecewise-lti}.
Conceptually, this system consists of a set of \(K\) linear subsystems as
given by \eqref{eq:continuous-lti}.
At any given time \(t\), the state is determined by subsystem
\(l(t)\).
Hollanders assumes that the system switches between subsystems
instantaneously.
Therefore, \(l(t)\) is constant in between two moments when the system
switches from one subsystem to another.

\begin{equation}
  \label{eq:continuous-piecewise-lti}
  \dot{x}(t) = A_{l(t)}x(t)+B_{l(t)}u(t)+\xi_t
\end{equation}

To find the way in which gene expression levels influence each other,
Hollanders uses a system identification approach to find a system of
either the form \eqref{eq:continuous-lti} or one of the form
\eqref{eq:continuous-piecewise-lti} that fits the data.
In order to find a unique solution, Hollanders' approach uses
a trade-off between minimizing an error criterion and criteria based
on norms of the resulting system matrices.

A common feature of the papers we described above is that each paper
proposes a method to find the structure of a dynamical system from
input/output data. 
By the structure of a dynamical system, we mean a description of which
variable directly influences another variable.
This structure can be represented as a directed graph, in which each
vertex represents a variable and the edge \((x,y)\) exists in the
graph if the variable represented by \(x\) directly influences the
variable represented by \(y\).

In the vector AR models of Goebel et al., we assume that the
covariance matrix of the noise term is constant.
This implies that the measures for Granger causality, which are based
on these matrices, are constant for a given set of models.
Therefore, the structure of a vector AR model as proposed by Goebel et
al.\ will be constant.

For a linear system given by \eqref{eq:continuous-lti}, one variable
influences another if and only if a certain entry of either \(A\) or
\(B\) is non-zero.
Since the matrices \(A\) and \(B\) are constant, the structure of a
linear system is also constant.
This implies that the structure of each subsystem of a piecewise
linear system is constant.
Thus, a piecewise linear system must have one of a given set of
structures at any moment.
This structure changes only when the system switches from one
subsystem to another.

In bilinear systems such as those used by Friston et al., the
structure of the system is determined by the input signal. 
Therefore, the structure can change from moment to moment, dictated by
changes in input.
The way in which the inputs dictate the structure is, however,
constant.

From the above discussion, we see that for a given system, either the
structure is constant from moment to moment or a description of the
possible structures exists.
For the identification procedures of Goebel et al.\ , Friston et
al.\ and Hollanders to be practically useful, the possible structures
a system can have should be uniquely determined by input/output data.
As mentioned, the goal of this paper is to describe to what degree
this is the case for a linear system.

\subsection{Related work}
Before we consider how we will proceed and state our main
conclusions,we note that related work has  been done in the area of linear structured systems.
Essentially, a linear structured system as used by Dion et al.\cite{dion-et-al} and
other researchers is identical to what we call the structure of a given linear system.
The question considered by researchers in this area is what
information about systems with a particular structure can be derived
from the structure itself.
As Dion et al.\ indicate, a number of properties will either not hold
at all for systems with a given structure or hold for almost all
systems with this structure.
A property that holds for most systems with a structure is said to be
generic for this system.
Dion et al.\ give graph-theoretical criteria for genericity of a
number of properties, such as controllability, observability and the
solvability of a number of other control problems.
Though this work is useful, it does not answer our basic question,
namely how a given input/output relation determines the structure of
systems realizing this relation.

\subsection{Contribution}
In this paper, we consider to what degree the structure of a linear
system can be determined using the input/output data it generates.
First, we consider how linear transformations affect the structure of a linear system.
Here, we see that for many linear systems, there exists a linear
transformation that changes the structure of the system.
Unless the linear transformation is part of a relatively small
set, we have no guarantee that it does not change the structure
of the system.

In light of this result, we consider whether a weaker
graph-theoretical relation than isomorphism could allow us to find
some aspect of the structure that is conserved by linear
transformations.
For this purpose, we consider a number of variants of graph
homomorphism.
Unfortunately, each variant is either not an equivalence relation or
leads to strange results where systems that should have different
structures have the same structure.

Given that linear transformations that result in isomorphic graphs
appear to be rare and variants of homomorphism lead to undesirable
results, we apply graph isomorphism to condensed graphs.
The condensed graph of a system is obtained from the original graph
of the system by replacing each strong component by a single vertex.
An edge exists from one component to another if a vertex in the
first component had an edge to a vertex in the second component in
the original graph.
Unfortunately, the same problem we saw with graph isomorphism
applied directly to a system's graph also occurs with
condensed graphs.
That is, there exist linear transformations that result in systems
with different condensed graphs.

The results discussed above indicate that if we do not make any
assumptions about the system, very few, if any, aspects of the
system's structure are determined by input/output data.
Thus, it seems that we need to make assumptions about a system in
order to determine its structure using this data.

The first kind of assumption we consider is that the system is a
minimal SISO system in some canonical form.
In a number of special cases, these assumptions allow us to
characterize the existence of isomorphisms between the graphs of two
systems. In other words, we can characterize when two systems have the
same structure.

The second kind of assumption we consider is that some variables do
not influence certain other variables.
That is, we assume that some edges do not occur in the graph of our
system.
This implies that the graph of our system must be a subgraph of a
given graph.
Using this graph, a number of properties of systems satisfying our
assumptions can be derived, as described by Dion et al.
Our first result concerning these properties is that unless there is
no edge from a state variable to an output in the given graph, we
cannot uniquely identify our system from input/output data.
This implies that for a realistic system, where the state variables do
influence the output, we cannot find the system parameters using
input/output data.
In our second result, we give a graph-theoretical characterization of
graphs such that almost all systems satisfying the assumptions given
by this graph are minimal.
These conditions are also necessary conditions for any given system to
be minimal.

\subsection{Application of the results}
Since the results described above apply only to linear systems, they
are not directly applicable to the models considered by Friston et
al.\cite{friston}, Goebel et al.\cite{goebel} and
Hollanders\cite{hollanders}.
However, the linear systems we consider are strongly related to each
of the model classes used by these researchers.
For instance, if we set each matrix \(B_i\) to zero in
\eqref{eq:bilinear-state}, we find the state equation of a linear
system.
Thus, we may identify the resulting system with a linear system where
the output of the system equals its state.
Due to the assumption that the system's output equals its state, we
cannot directly apply a linear transformation to such a system.
If we instead allow the output to be any linear
transformation of the state, we find a subset of the bilinear systems
corresponding to the linear systems such that the input does not
directly influence the output.
Many of our results apply to this set of linear systems.
Since this set of systems is a proper subset of the set of all
bilinear systems, we conjecture that these results may apply, possibly
in a modified form, to the class of all bilinear systems.

A similar relation exists between the linear systems we consider and
the piecewise linear systems used by Hollanders.
Indeed, if we ignore the noise term in \eqref{eq:continuous-lti} and
allow the output of the system to be a linear combination of its
state, we find the same set of linear systems we considered above.
In this case, we can generalize to piecewise linear systems by
allowing the system to have an arbitrary number of subsystems.
Since the systems with 1 subsystem are a subset of this more
general class, we conjecture that many of our results may also hold
for this more general case.

The relation between the vector AR models considered by Goebel et
al.\cite{goebel} and our linear systems is more complex.
Though we can represent each vector AR model by a linear system driven
by a stochastic input, our notion of this system's graph structure
may differ from the Granger causality criteria used by Goebel et al.
In some cases, these structures can give equivalent results.
For instance, suppose we have two time series \(x_i \in \mathbb{R}^m\)
and \(y_i \in \mathbb{R}^n\).
Suppose the best vector AR model of \(x_i\) and \(y_i\) of the form
\eqref{eq:vector-ar} has block-diagonal matrices \(A_i\), where \(A_i
= 
\begin{bmatrix}
  A_{i,x} & 0 \\
  0 & A_{i,y}
\end{bmatrix}\), \(A_{i,x} \in \mathbb{R}^{m \times m}\) and \(A_{i,y}
\in \mathbb{R}^{n \times n}\),
and the
covariance matrix of the noise \(u_n\)  is also block-diagonal,
i.e.\ this matrix has the form
\(\begin{bmatrix}
  \Sigma_x & 0 \\
  0 & \Sigma_y 
\end{bmatrix}\), where \(\Sigma_x \in \mathbb{R}^{m \times m}\) and
\(\Sigma_y \in \mathbb{R}^{n \times n}\).
In this case, the graph of this system consists of two disjoint
subgraphs, one for the vector AR model of \(x_i\) and one for the
model of \(y_i\).
Since the joint model of \(x_i\) and \(y_i\) predicts each of these
time series using only its own past values, we also find that there is
no Granger causality between \(x_i\) and \(y_i\).
Thus, we see that the graph structure and Granger causality criteria
can give equivalent results.
It is unclear at this point whether and, if so, how, this holds in
the general case.


The remainder of this paper is structured as follows.
In Section \ref{sec:systems-graph}, we recall the definitions of a
linear system, its associated graph and other relevant concepts.
We then use these concepts in Section \ref{sec:main-results} to state
the results we have discussed above.
These results are proved in Section \ref{sec:proofs-main-results}.
In the final section, we briefly review our main results and their
implications.

\section{Systems, graph structure and equivalent structures}
\label{sec:systems-graph}
In this section, we recall the definition of a discrete-time LTI
system and its associated graph structure.
Given an initial state \(x_0\) and an input sequence \(u_k\), the
evolution over time of the state \(x_k\) and output \(y_k\) of the
system is given by \eqref{eq:system-state} and \eqref{eq:system-output}.
In these equations, the matrices \(A \in \mathbb{R}^{n_x \times
  n_x}\), \(B \in \mathbb{R}^{n_x \times n_u}\), \(C \in
\mathbb{R}^{n_y \times n_x}\) and \(D \in \mathbb{R}^{n_y \times n_u}\)
 are  constants.
Since these constants uniquely determine the discrete-time LTI system,
we define such a system to be a 4-tuple of the matrices
\(A\),\(B\), \(C\) and \(D\), as in Definition \ref{def:dt-lti} below.
Similar definitions can be found in many textbooks on linear systems,
such as those by Vaccaro \cite{vaccaro} and Kailath \cite{kailath}.

\begin{align}
\label{eq:system-state} x_{k+1} = Ax_k+Bu_k \\
\label{eq:system-output} y_k = Cx_k + Du_k
\end{align}

\begin{definition}
\label{def:dt-lti}
An LTI state-space system with \(n_u\) inputs, \(n_y\) outputs and
\(n_x\) state variables is a 4-tuple of matrices
\(\left(A,B,C,D\right)\), where \(A \in \mathbb{R}^{n_x \times n_x}, B \in \mathbb{R}^{n_x \times n_u}, C \in \mathbb{R}^{n_y \times n_x}, D \in \mathbb{R}^{n_y \times n_u}\)
\end{definition}

We define the structure of a state-space system as a directed graph
whose vertices are the inputs, outputs and state variables of the
system. 
An edge from one variable \(x\) to another variable \(y\) exists in
this graph if \(y\) is directly influenced by \(x\). 
We formalize this concept below.
We begin by defining the graph of a matrix and then use this concept
to define the graph of a system.
Before we state these definitions, we briefly recall the concept of a
directed graph as defined in many textbooks on graph theory, such as
the book by Chartrand et al.\cite[Ch.~ 7]{chartrand}.

\begin{definition}
  \label{def:digraph}
  A (directed) graph \(G\) is a 2-tuple \((V,E)\), where \(V\) is a
  finite set of objects called vertices and \(E \subset V \times
  V\) is a set of 2-tuples of vertices called edges.
\end{definition}

\begin{definition}
  \label{def:associatedgraph-matrix}
  Let \(M \in \mathbb{R}^{n \times n}\).
  Then, the associated graph \(G\) of the matrix \(M\) is given by \(G
  = (V,E)\), where \(V = \{v_1,v_2,\cdots,v_n\}\) and \(E =
  \{(v_i,v_j) | M_{ji} \not = 0\}\).
\end{definition}

\begin{definition}
  \label{def:associated-graph}
  The associated graph of the system \stdnamedsys with \(n_u\) inputs,
  \(n_x\) state variables and \(n_y\) outputs is defined as the
  graph \(G(S) = (V,E) = \left(\{v_1,v_2,\cdots,v_n\},E\right)\) of the matrix \(M_s = 
\begin{bmatrix}
A & B & 0 \\
0 & 0 & 0 \\
C & D & 0
\end{bmatrix}\), where \(n = n_x + n_u + n_y\) and \(M_s \in
\mathbb{R}^{n \times n}\).
  Additionally, we will use the following notation.
  The vertices \(v_1\) through \(v_{n_x}\) are called the state
  variable vertices of \(G(S)\).
  We will denote these vertices by \(x_i(G(S)) = v_i\).
  The vertices \(v_{n_x+1}\) through \(v_{n_x+n_u}\) are called the
  input vertices of \(G(S)\), denoted by \(u_i(G(S)) = v_{n_x+i}\).
  Finally, the remaining vertices are called the output vertices
  \(y_i(G(S))\), given by \(y_i(G(S)) = v_{n_x+n_u+i}\).
\end{definition}
\begin{remark}
  Unless otherwise noted, \(G(S)\) denotes the graph of a linear
  system as defined in Definition \ref{def:dt-lti}.
  When the graph we are referring to is clear from the context, we
  will abbreviate \(x_i(G(S))\) to \(x_i\).
  We will similarly abbreviate \(y_i(G(S))\) to \(y_i\) and \(u_i(G(S))\) to \(u_i\).
\end{remark}

A useful corollary of Definition \ref{def:associated-graph} is that
only certain edges occur in the graph \(G(S)\). Furthermore, these
edges occur if and only if the corresponding entries in the system
matrices of \(S\) are non-zero.
We state this formally below.

\begin{corollary}
\label{cor:associated-graph}
Let \stdnamedsys be a linear system with \(n_u\) inputs, \(n_x\) state
variables and \(n_y\) outputs.
Then, in the graph \(G(S) = (V,E)\) of \(S\), each edge \((v_1,v_2)\) is of
one of the forms listed below.
\begin{enumerate}
\item  \((x_j,x_i)\),where \(1 \leq i,j \leq n_x\)
\item \((u_j,x_i)\), where \(1 \leq i \leq n_x \wedge 1 \leq j \leq n_u\)
\item   \((x_j,y_i)\), where \(1 \leq i \leq n_y \wedge 1 \leq j \leq n_x\) 
\item  \((u_j,y_i)\), where \(1 \leq i \leq n_y \wedge 1 \leq j \leq n_u\)
\end{enumerate}
In addition, the following conditions hold.
\begin{enumerate}
\item \(\left(x_j,x_i\right) \in E\) if and only if \(A_{ij} \not = 0\)
\item \(\left(u_j,x_i\right) \in E\) if and only if \(B_{ij} \not = 0\)
\item \(\left(x_j,y_i\right) \in E\) if and only if \(C_{ij} \not = 0\)
\item \(\left(u_j,y_i\right) \in E\) if and only if \(D_{ij} \not = 0\)
\end{enumerate}
\end{corollary}

In this paper, we will consider a number of functions from the vertex
set of one graph to the vertex set of another.
We will restrict these functions to only map vertices onto vertices ``of the same kind''. 
That is, inputs must be mapped onto inputs, outputs onto outputs and so on.
Formally, we say that the function must be type-restricted, as defined in Definition \ref{def:type-restricted} below.


\begin{definition}
Let \(S_1\) be a linear system with \(n_{x,1}\) state variables,
\(n_{u,1}\) inputs and \(n_{y,1}\) outputs and \(S_2\) a linear system
with \(n_{x,2}\) state variables, \(n_{u,2}\) inputs and \(n_{y,2}\) outputs.
Then, the vertices \(v_1 \in V(G(S_1))\) and \(v_2 \in V(G(S_2))\) are
of the same type if one of the following conditions holds for some
\(i\) and \(j\):
\begin{enumerate}
\item \(1 \leq i \leq n_{u,1} \wedge 1 \leq j \leq n_{u,2} \wedge u_i(G(S_1)) = v_1 \wedge u_j(G(S_2)) = v_2\)
\item \(1 \leq i \leq n_{y,1} \wedge 1 \leq j \leq n_{y,2} \wedge y_i(G(S_1)) = v_1 \wedge y_j(G(S_2)) = v_2\)
\item \(1 \leq i \leq n_{x,1} \wedge 1 \leq j \leq n_{x,2} \wedge x_i(G(S_1)) = v_1 \wedge x_j(G(S_2)) = v_2\)
\end{enumerate}
\end{definition}

\begin{definition}
\label{def:type-restricted}
Let \(S_1\) and \(S_2\) be linear systems and let \(G(S_i) =
(V_i,E_i)\) for \(i = 1,2\).
A function \(\phi : V_1 \rightarrow V_2\) is said to be type-restricted if, for all \(v \in V_1\), \(v\) and \(\phi(v)\) are of the same type.
\end{definition}

We will now illustrate these definitions with an example.

\begin{figure}
\caption{The associated graph \(G(S)\) in Example \ref{ex:graphs-and-mappings}}
\label{fig:graph-example}
\begin{center}
\includegraphics{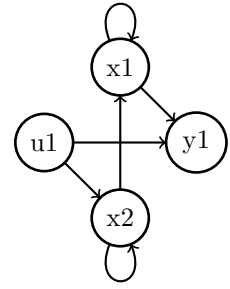}
\end{center}
\end{figure}

\begin{example}
\label{ex:graphs-and-mappings}
Consider the system \(S = \left(
\begin{bmatrix}
1 & 2 \\
0 & 1 
\end{bmatrix},
\begin{bmatrix}
0 \\
3
\end{bmatrix}
,
\begin{bmatrix}
1 & 0
\end{bmatrix},
\begin{bmatrix}
2  
\end{bmatrix}
\right)\).
The associated graph \(G(S)\) of \(S\) is shown in Figure \ref{fig:graph-example}.
Of the two functions \(\phi_1\) and \(\phi_2\) defined below,
\(\phi_1\) is type-restricted while \(\phi_2\) is not.
\begin{align*}
\phi_1(v) & = \begin{cases}
x_2(G(S)) & \text{if } v = x_1(G(S)) \\
x_1(G(S)) & \text{if } v = x_2(G(S)) \\
v & \text{otherwise}
\end{cases} \\
\phi_2(v) & = \begin{cases}
u_1(G(S)) & \text{if } v = x_1(G(S)) \\
x_1(G(S)) & \text{if } v = u_1(G(S)) \\
v & \text{otherwise}
\end{cases}
\end{align*}
\end{example}

The first kind of function we will consider is an isomorphism.
A (directed-)graph isomorphism is defined in Definition \ref{def:graph-isomorphism}.
Unless otherwise noted, we will additionally require that the
isomorphism is type-restricted as defined above.

\begin{definition}
\label{def:graph-isomorphism}
Let \(S_1\) and \(S_2\) be linear systems and let \(G(S_i) =
(V_i,E_i)\) for \(i = 1,2\).
A type-restricted function \(\phi : V_1 \rightarrow V_2\)  is a type-restricted isomorphism if:
\begin{enumerate}
\item \(\phi\) is a bijection, that is, it is both injective and
  surjective.
\item For any vertices \(u,v \in V_1\), \((u,v) \in E_1\) if and
  only if \((\phi(u),\phi(v)) \in E_2\).
\end{enumerate}
\end{definition}

A type-restricted isomorphism as defined above may permute the inputs
and outputs of a system.
The results we derive in the remainder of this paper remain valid if
we require a type-restricted isomorphism to leave the inputs and
outputs in the same order.
Thus, if desired, this condition may be added to Definition
\ref{def:graph-isomorphism}.

In addition to the isomorphisms we have defined above, we will find it
useful to consider isomorphisms on condensed graphs.
The condensed graph corresponding to a given graph is the graph whose
vertices are the components of the original graph.
An edge from one component to another exists in the condensed graph if
an edge from a vertex in one component to a vertex in another existed
in the original graph.
To formally state the definition of a condensed graph and a
condensed-graph isomorphism, we begin by stating the definition of a
component.

\begin{definition}
  Let \(S\) be a linear system and let \(G(S) = (V,E)\) be its
  associated graph.
  For all \(a,b \in V\), let \(a \leftrightarrow b\) if and only if
  there exist directed paths from \(a\) to \(b\) and vice-versa in \(G(S)\).
  The relation \(\leftrightarrow\) defined above is an equivalence
  relation.
  The equivalence classes of this relation are called the strong components
  of \(G(S)\).
\end{definition}

From the definition above, it is clear that the strong components of
\(G(S) = (V,E)\) are sets of vertices that form a partition of \(V\).
In order to define the condensed graph \(CG(S)\) of \(S\) and
type-restricted isomorphisms on such graphs, we will need to assign a
type to each component of \(G(S)\).
To do this, we recall an observation previously made by
Dion et al.\ in their survey of structured linear systems
\cite{dion-et-al}.
A formal proof of this observation is given in Appendix \ref{app:prelims}.

\begin{observation}
  \label{obs:input-output-singletons}
  Each input vertex \(u_i(G(S))\) or output vertex \(y_i(G(S))\) is
  the only element of its component in \(CG(S)\).
\end{observation}

By the above observation, any component of \(G(S)\) with two or more
elements must consist entirely of state variables.
We will use this fact below to state the definition of a condensed
graph and type-restricted mappings between such graphs.

\begin{definition}
Let \(S\) be a linear system with \(n_u\) inputs and \(n_y\) outputs
and let \(G(S) = \left(V_G,E_G\right)\) be its associated graph.
Furthermore, let \(c_1,c_2,\cdots,c_m\) be the strong components of \(G(S)\) with
more than one member.
Then, we define the vertex set \(V_{CG}\) by \(V_{CG} =
\{u_1(G(S)),u_2(G(S)),\cdots,u_{n_u}(G(S)),\allowbreak c_1,\allowbreak
c_2,\allowbreak\cdots,c_m,\allowbreak y_1(G(S)),y_2(G(S)),\cdots,y_{n_y}(G(S))\}\).
The edge set \(E_{CG}\) contains the elements specified by the conditions
below.
The condensed graph \(CG(S)\) is the graph \(CG(S) =
(V_{CG},E_{CG})\).

\begin{enumerate}
\item For all \(1 \leq i \leq n_u\) and \(1 \leq j \leq m\), \((u_i,c_j) \in E_{CG}\) if and only if a vertex \(v \in c_j\)
  exists such that \((u_i,v) \in E_G\).
\item For all \(1 \leq i,j \leq m\), \((c_i,c_j) \in E_{CG}\) if and only if vertices \(v_1 \in
  c_i\) and \(v_2 \in c_j\) exist such that \((v_1,v_2) \in E_G\).
\item For all \(1 \leq i \leq m\) and \(1 \leq j \leq n_y\), \((c_i,y_j) \in E_{CG}\) if and only if a vertex \(v \in
  c_i\) exists such that \((v,y_j) \in E_G\).
\item For all \(1 \leq i \leq n_u\) and \(1 \leq j \leq n_y\), \((u_i,y_j) \in E_{CG}\) if and only if \((u_i,y_j) \in E_G\).
\end{enumerate}

In addition, we will denote the vertices \(u_i(G(S))\) by
\(u_i(CG(S))\) and call these vertices input vertices.
Similarly, the vertices \(y_i(G(S))\) will be denoted by
\(y_i(CG(S))\) and will be called output vertices.
Finally, we will call the vertices \(c_i\) state variable component
vertices and denote them by \(c_i(CG(S))\).
\end{definition}

\begin{definition}
  Let \(S_1\) be a linear system with \(n_{u,1}\) inputs and \(n_{y,1}\)
  outputs and \(S_2\) a linear system with \(n_{u,2}\) inputs
  and \(n_{y,2}\) outputs.
  Furthermore, let \(n_{c,i}\), \(i = 1,2\), be the number of state
  variable component vertices of \(CG(S_i)\).
  Finally, let \(CG(S_i) = (V_i,E_i)\), for \(i = 1,2\).
  Then, a mapping \(\phi : V_1 \rightarrow V_2 \) is
  type-restricted if, for all \(v \in V_1\), one of the following
  conditions holds for some \(i\) and \(j\).
  \begin{enumerate}
  \item \(1 \leq i \leq n_{u,1} \wedge 1 \leq j \leq n_{u,2} \wedge v = u_i(CG(S)) \wedge \phi(v) = u_j(CG(S'))\)
  \item \(1 \leq i \leq n_{y,1} \wedge 1 \leq j \leq n_{y,2} \wedge v = y_i(CG(S)) \wedge \phi(v) = y_j(CG(S'))\)
  \item \(1 \leq i \leq n_{c,1} \wedge 1 \leq j \leq n_{c,2} \wedge v = c_i(CG(S)) \wedge \phi(v) = c_j(CG(S'))\)
  \end{enumerate}
\end{definition}

Using the above definitions, we can define isomorphisms on condensed
graphs, which we will refer to as CG-isomorphisms.
As we define below, two systems \(S\) and \(S'\) such that a
CG-isomorphism exists between their condensed graphs are called
CG-isomorphic.

\begin{definition}
  Let \(S\) and \(S'\) be linear systems.
  Then, a type-restricted isomorphism between \(CG(S)\) and \(CG(S')\) is
  said to be a condensed-graph (CG-) isomorphism between \(S\) and
  \(S'\).
  If a CG-isomorphism between \(S\) and \(S'\) exists, \(S\) and
  \(S'\) are said to be CG-isomorphic, written \(S \simeq_{CG} S'\).
\end{definition}

\section{Main results}
\label{sec:main-results}
\subsection{Graph isomorphism and its inadequacy}
\label{subsec:main-graph-isomorphism}
We begin by defining linear transformations of systems and considering
how these transformations affect the associated graph of the system.
Consider the system \namedsys{S}{A}{B}{C}{D} with \(n_u\) inputs,
\(n_x\) state variables and \(n_y\) outputs.
The evolution over time of the state vector \(x_k \in
\mathbb{R}^{n_x}\) and output \(y_k \in \mathbb{R}^{n_y}\)
of \(S\), given an input \(u_k \in \mathbb{R}^{n_u}\), is given by \eqref{eq:system-state} and \eqref{eq:system-output}.
Let \(T\) be an invertible \(n_x \times n_x\) matrix and \(z_k = Tx_k\).
Then, the evolution over time of \(z_k\) and \(y_k\) is given by 
\eqref{eq:trans-system-state} and \eqref{eq:trans-system-output}.
We define the result of transforming the system \(S\) using the matrix
\(T\) to be the system with state vector \(z_k\).
This definition is stated formally below.
\begin{align}
\label{eq:trans-system-state} z_{k+1} & = TAT^{-1}z_k+TBu_k\\
\label{eq:trans-system-output} y_k & = CT^{-1}z_k + Du_k
\end{align}

\begin{definition}
The result of transforming the system \namedsys{S}{A}{B}{C}{D} with
\(n_x\) state variables using an invertible \(n_x \times n_x\) matrix \(T\), denoted \transys{S}{T}, is 
the system \((TAT^{-1},TB,CT^{-1},D)\).
\end{definition}

Our first result is that in many cases, we can find a new system with
a non-isomorphic graph by diagonalizing the \(A\)-matrix of the
system.
Thus, it is possible for a linear transformation to result in a system
with a different graph.

\begin{theorem}
\label{thm:diagonalizable-non-isomorphism}
For a system \stdnamedsys, where \(A\) is diagonalizable but not
diagonal, there exists an invertible matrix \(T\) such that \(G(S)\)
and \(G(\transys{S}{T})\) are non-isomorphic.
\end{theorem}

On the other hand, there do exist matrices \(T\) such that
\(G(\transys{S}{T})\) and \(G(S)\) are isomorphic for all \(S\).
For any positive integer \(n\), we call the set of all such \(n \times n\) matrices \(GI(n)\).
The set \(GI(n)\) is a subgroup of the group of invertible \(n \times
n\) matrices, as stated below.

\begin{definition}
\label{def:gi-n}
The set \(GI(n)\) consists of the \(n \times n\) invertible matrices
\(T\) such that \(G(S)\) and \(G(\transys{S}{T})\) are isomorphic for
all linear systems \(S\) with \(n\) state variables.
\end{definition}

\begin{theorem}
\label{thm:gi-n-subgroup}
\(GI(n)\) is a subgroup of all invertible \(n \times n\) matrices.
Equivalently, the following conditions hold:
\begin{enumerate}
\item \(I_n \in GI(n)\)
\item \(GI(n)\) is closed under matrix multiplication
\item \(GI(n)\) is closed under matrix inversion.
\end{enumerate}
\end{theorem}

Though the above result shows that \(GI(n)\) has some algebraic
structure, it does not indicate how many and which elements this group
has.
Below, we state two results that show that certain classes of matrices
are subsets of \(GI(n)\).
One of these classes, the permutation matrices, is defined below.
In Theorem \ref{thm:diag-gi-n}, \(\{m_i\}\) denotes the diagonal
matrix whose diagonal elements are \(m_i\).
Thus, \(\{m_1,m_2,\cdots,m_n\}\) denotes the matrix \(M\) given below.

\[ M = 
\begin{bmatrix}
  m_1 & 0 & \cdots & 0 \\
  0 & m_2 & \cdots & 0 \\
  \vdots & \vdots & \ddots & \vdots \\
  0 & 0 & \cdots & m_n
\end{bmatrix}\]

\begin{theorem}
\label{thm:diag-gi-n}
For any diagonal matrix \(M = \{m_i\}\), \( 1 \leq i
\leq n\), such that \(m_i \not = 0\) for all \(i\), \(M \in GI(n)\).
\end{theorem}

\begin{definition}
For a permutation \(e_1,e_2,\cdots,e_n\) of the integers
\(1,2,\cdots,n\), the corresponding \(n \times n\) permutation matrix
is the matrix \(P(e_1,e_2,\cdots,e_n) = (P_{ij})\), \(1 \leq i,j \leq n\), where \(P_{ij}\) is
as given below.
\[ P_{ij} = 
\begin{cases}
1 & \text{if } j = e_i \\
0 & \text{otherwise}
\end{cases}
\]
\end{definition}

\begin{theorem}
\label{thm:permutation-gi-n}
Any \(n \times n\) permutation matrix is an element of \(GI(n)\).
\end{theorem}
\subsection{Graph homomorphism}
In most of this paper, we concentrate on graph isomorphism.
In principle, we could also have focused on a weaker relation, such as
graph homomorphism.
In this subsection, we state a number of results that indicate why we
consider graph homomorphism to be inadequate for our purposes.
We begin by defining a graph homomorphism below.
Intuitively, a graph homomorphism is a mapping from one graph to
another that conserves edges.

\begin{definition}
  Let \(S_1\) and \(S_2\) be linear systems.
  Furthermore, let \(G(S_i) = (V_i,E_i)\) for \(i = 1,2\).
  Then, a type-restricted function \(\phi : V_1 \rightarrow V_2\) is a
  type-restricted homomorphism if for all \(u,v \in V_1\), \((u,v) \in
  E_1\) implies \((\phi(u),\phi(v)) \in E_2\).
\end{definition}

Our first example shows that the graph of every system is homomorphic
to that of a standard system.

\begin{figure}
  
  \caption{The system \(S_1\) from Example \ref{ex:homomorphism}}
  \label{fig:homomorphism-ex}
  \begin{center}
    \includegraphics{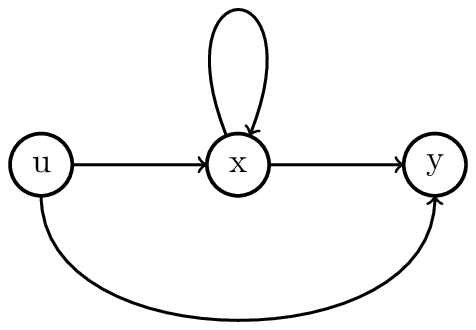}
  \end{center}
\end{figure}

\begin{example}
 \label{ex:homomorphism}
 Consider the system \(S_1 = (
 \begin{bmatrix}
   1
 \end{bmatrix},
 \begin{bmatrix}
   1
 \end{bmatrix},
 \begin{bmatrix}
   1
 \end{bmatrix},
 \begin{bmatrix}
   1
 \end{bmatrix})\).
 The graph of \(S_1\) is shown in Figure \ref{fig:homomorphism-ex}.
 Next, consider an arbitrary system \stdnamedsys.
 The mapping \(\phi\) defined below is a type-restricted homomorphism
 from \(S\) to \(S_1\).
 \[
 \phi(v) = \begin{cases}
   u_1(G(S_1)) & \text{ if } v = u_i(G(S)) \text{ for some } i \\
   y_1(G(S_1)) & \text{ if } v = y_i(G(S)) \text{ for some } i \\
   x_1(G(S_1)) & \text{ if } v = x_i(G(S)) \text{ for some } i
   \end{cases}
 \]
\end{example}

If graph homomorphism were an equivalence relation, the above example
would be problematic since every system would have the same structure.
Fortunately, graph homomorphism fails to be symmetric, and so cannot
be an equivalence relation.
This fact is shown in the following example.

\begin{figure}
  
  \caption{The system \(S_2\) from Example \ref{ex:homomorphism}}
  \label{fig:homomorphism-ex2}
  \begin{center}
    \includegraphics{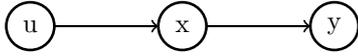}
  \end{center}
\end{figure}

\begin{example}
  We will now show that, though a homomorphism from \(S \) to \(S_1\)
  exists for all \(S\), the opposite is not necessarily the case.
  To do this, we consider the system \(S_2 = (
  \begin{bmatrix}
    0
  \end{bmatrix},
  \begin{bmatrix}
    1
  \end{bmatrix},
  \begin{bmatrix}
    1
  \end{bmatrix},
  \begin{bmatrix}
    0
  \end{bmatrix})\), whose graph is shown in Figure
  \ref{fig:homomorphism-ex2}.
  Since any homomorphism \(\phi\) from \(S_1\) to \(S_2\) must be
  type-restricted, it follows that \(\phi(x_1(G(S_1))) =
  x_1(G(S_2))\).
  But then, by the definition of a type-restricted homomorphism,
  \((x_1(G(S_2)),x_1(G(S_2)))\) should be an edge of \(G(S_2)\), since
  \((x_1(G(S_1)),x_1(G(S_1)))\) is an edge of \(G(S_1)\).
  Since this is not the case, no type-restricted homomorphism from
  \(S_1\) to \(S_2\) can exist.
\end{example}

One way to make homomorphism an equivalence relation is to require not
just that a homomorphism from \(G_1\) to \(G_2\) exists, but also that
a homomorphism from \(G_2\) to \(G_1\) exists.
Even with this interpretation, homomorphism can lead to strange
results, as we show below.

\begin{figure}
  
  \caption{The system \(S_3\) from Example
    \ref{ex:homomorphism-structure}.}
  \begin{center}
    \includegraphics{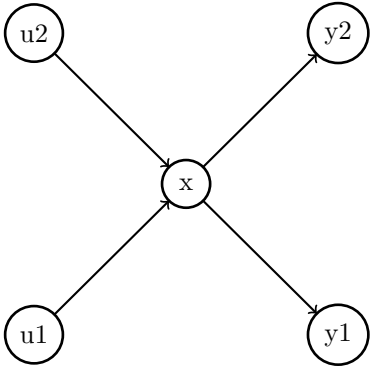}
  \end{center}
\end{figure}

\begin{example}
  \label{ex:homomorphism-structure}
  Our next example concerns the system \namedsys{S_3}{
    \begin{bmatrix}
      0
    \end{bmatrix}}
  {
    \begin{bmatrix}
      1  & 1\\
    \end{bmatrix}
  }
  { 
    \begin{bmatrix}
      1 \\
      1
    \end{bmatrix}
  }
  {
    \begin{bmatrix}
      0
    \end{bmatrix}
  }.
The mapping \(\phi\) below is a homomorphism from \(S_3\) to \(S_2\),
and similarly, \(\psi\) is a homomorphism in the other direction.
\begin{align*}
\phi(v) &= \begin{cases}
u_1(G(S_2)) & \text{ if } v = u_1(G(S_3)) \vee v = u_2(G(S_3)) \\
y_1(G(S_2)) & \text{ if } v = y_1(G(S_3)) \vee v = y_2(G(S_3)) \\
x_1(G(S_2)) & \text{ if } v = x_1(G(S_3))
\end{cases} \\
\psi(v) &= \begin{cases}
u_1(G(S_3)) & \text{ if } v = u_1(G(S_2)) \\
y_1(G(S_3)) & \text{ if } v = y_1(G(S_2)) \\
x_1(G(S_3)) & \text{ if } v = x_1(G(S_2))
\end{cases}
\end{align*}
So, if we consider two systems to have the same structure if
homomorphisms exist between them in both directions, \(S_2\) and
\(S_3\) have the same structure.
\end{example}

In the above example, we show that the systems \(S_2\) and \(S_3\)
have the same structure.
Intuitively, this should not be the case, since \(S_2\) and \(S_3\)
have different numbers of inputs.
From the above examples, we conclude that homomorphism, at least in
the variants we have discussed here, is inadequate for our purposes
since it will always lead to strange results.

\subsection{Condensed-graph isomorphism}
\label{subsec:main-cg-isomorphism}
We will now consider how some of the results we derived earlier for
graph isomorphism apply to condensed-graph isomorphism.
Our first result is that the counterexample we presented for graph
isomorphism also applies to condensed-graph isomorphism.

\begin{theorem}
  \label{thm:condensed-diag-counterexample}
  A system \stdnamedsys with \(A\) non-diagonal but diagonalizable, is
  not CG-isomorphic to \transys{S}{T}, where \(D = TAT^{-1}\) is diagonal.
\end{theorem}

In our next result, we show that there exists a class of systems whose
condensed graphs contain 1 state variable component.
Each of these systems has a diagonalizable \(A\)-matrix, allowing us
to transform each system to another system where each state variable is in a
component of its own.

\begin{theorem}
  \label{thm:1-n-components}
  For each integer \(n\), there exists a system \stdnamedsys with
  \(n\) state variables such that its condensed graph \(CG(S)\)
  contains 1 state variable component.
  In addition, the matrix \(A\) is diagonalizable but not diagonal,
  implying that there exists a system \transys{S}{T} such that
  \(CG(\transys{S}{T})\) has \(n\) state variable components.
\end{theorem}

The results above indicate that the condensed-graph structure of a
system is not conserved by linear transformations.
Since linear transformations do not change the input/output behavior
of a system, this implies that systems with different condensed-graph
structures may produce the same input/output data.
Therefore, we cannot identify the condensed-graph structure of a
system using input/output data unless we have additional
assumptions about the system.
In other words, to determine the (condensed-)graph structure of a
system, we need a parameterization of the set of linear systems we are
interested in.
Given a particular parameterization, we can ask when two systems in
this parameterization have isomorphic structures.
We will now consider this question for two parameterizations.
The first parameterization concerns the minimal SISO systems with
diagonal \(A\)-matrices.
In the second parameterization, we consider minimal SISO systems whose
\(A\)-matrices are in the natural normal forms given by Gantmacher \cite{gantmacher}.

For systems in the first parameterization, i.e. minimal SISO systems
with diagonal \(A\)-matrices, we can characterize the existence of a
type-restricted isomorphism between two systems.
This result is stated below.

\begin{theorem}
\label{thm:diag-isomorphism}
The graphs of two minimal diagonal SISO realizations
\stdnamedsys  and \namedsys{S'}{A'}{B'}{C'}{D'} will be
isomorphic if and only if the following conditions hold:
\begin{enumerate}
\item Either \(D = D' = 0\) or both \(D\) and \(D'\) are non-zero
\item \(A\) and \(A'\) have the same number of non-zero elements along
  their diagonals
\end{enumerate}
\end{theorem}

In the second parameterization, a system's \(A\)-matrix is in one of
the normal forms given by Gantmacher \cite{gantmacher}.
These normal forms are based on the elementary divisors and invariant
polynomials of the matrix \(A\) and are defined below.
The concepts of invariant polynomials and elementary divisors are
recalled in the Appendix.

\begin{definition}
Let \(l(s) = s^n + a_1s^{n-1}+a_2s^{n-2}+\cdots+a_n\).
The companion matrix \(L\) for the polynomial \(l(s)\) is the \(n
\times n\) square matrix shown below.
As shown by Gantmacher \cite[Ch.~6]{gantmacher}, \(\dcharmat{L} =
l(\lambda)\). Furthermore, all invariant polynomials \(i_j\) of \(L\) other
than \(i_1\) are equal to 1, and so \(i_1 = \dcharmat{L}\).
\[ L = \left[
\begin{array}{ccccc}
0 & 0 & \cdots & 0 & -a_n \\
1 & 0 & \cdots & 0 & -a_{n-1} \\
0 & 1 & \cdots & 0 & -a_{n-2} \\
\cdots & \cdots & \cdots & \cdots & \cdots\\
0 & 0 & \cdots & 1 & -a_1 
\end{array}
\right]\]  
\end{definition}

\begin{definition}
Let \(A\) be a real matrix and
\(i_1,i_2,\cdots,i_r,i_{r+1},\cdots\,i_n\) be its invariant
polynomials, such that the polynomials \(i_1\) through \(i_r\) are of
positive degree and \(i_j = 1\) for \(r+1\leq j \leq n\).
Then, the first natural normal form of A is the matrix \(L =
\left\{L_1,L_2,\cdots,L_r\right\}\)
\footnote{Here, \(\left\{A_1,A_2,\cdots,A_n\right\}\) denotes the
  block-diagonal matrix whose diagonal blocks are
  \(A_1,A_2,\cdots,A_n\).}, where \(L_j\) is the companion matrix for
the polynomial \(i_j\).
\end{definition}

\begin{definition}
Let \(e_1,e_2,\cdots,e_k\) be the elementary divisors of \(A\).
The second natural normal form of A is the matrix \(L =
\left\{L_1,L_2,\cdots,L_k\right\}\), where \(L_i\) is the companion matrix of
\(e_i\).
\end{definition}

\begin{remark}
  As Gantmacher \cite{gantmacher} shows, every matrix \(A\) is similar
  to a matrix \(A_1\) in first natural normal form and a matrix
  \(A_2\) in second natural normal form.
  Thus, for every linear system \stdnamedsys, there exist systems
  \namedsys{S_1}{A_1}{B_1}{C_1}{D_1} and
  \namedsys{S_2}{A_2}{B_2}{C_2}{D_2} such that \(A_1\) is in first
  natural normal form, \(A_2\) is in second natural normal form and
  both \(S_1\) and \(S_2\) are similar to \(S\).
\end{remark}

As we state below, if a SISO realization is minimal and its
\(A\)-matrix is in first natural normal form, the matrix \(A\) is a
companion matrix.
This implies that the system will be similar to the well-known
canonical forms of such systems.

\begin{theorem}
\label{thm:first-natural-companion}
If a SISO realization \((A,B,C,D)\) is minimal and the matrix \(A\) is in first natural normal form, \(A\) is a companion matrix.
\end{theorem}

Unlike first natural normal form, second natural normal form differs
from the well-known canonical forms.
In this form, the matrix \(A\) is block-diagonal and each block
corresponds to an elementary divisor of the matrix \(A\).
Below, we give a characterization of \(CG\)-isomorphism between two
realizations whose \(A\)-matrices are in second natural normal form.

\begin{theorem}
\label{thm:second-natural-isomorphism}
Let \namedsys{S_1}{A_1}{B_1}{C_1}{D_1} and
\namedsys{S_2}{A_2}{B_2}{C_2}{D_2} be two minimal SISO realizations
where \(A_1\) and \(A_2\) are in second natural normal form  
such that 0 is not an eigenvalue of either \(A_1\) or \(A_2\).
Then, \(G(S_1)\) and \(G(S_2)\) are CG-isomorphic if and only if:
\begin{enumerate}
\item The number of distinct irreducible polynomials that divide
  \dcharmat{A_1} equals the number of distinct irreducible polynomials
  that divide \dcharmat{A_2}
\item Either \(D_1 = D_2 = 0\) or \(D_1 \not = 0 \wedge D_2 \not = 0\)
\end{enumerate}
\end{theorem}

\subsection{Components of condensed graphs}
\label{subsec:main-cg-order}
In the previous section, we saw that the condensed-graph structure of
a system is not uniquely determined by its input/output behavior.
Thus, there may be a variety of systems \(S\) and associated graphs
\(G(S)\) corresponding to any set of input/output data.
In each of these graphs \(G(S)\), some sets of state variables do not
interact with any variables outside the set itself.
Such sets correspond to isolated components in the graph \(G(S)\).
By permuting the state variables of \(S\), we can ensure that the
variables in each component form a sequence of consecutive variables,
that is, they are the variables \(x_i,x_{i+1},\cdots,x_{i+j}\), for
some \(i\) and \(j\).
The resulting system will then have an \(A\)-matrix with a
block-diagonal structure, where each isolated component in the graph
corresponds to a diagonal block.
Conversely, if \(S\) has a block-diagonal \(A\)-matrix, each block of
\(A\) corresponds to an isolated component in \(G(S)\).
Thus, if we can find bounds on the number of diagonal blocks a
block-diagonal realization similar to \(S\) may have, the same bounds
should apply to the number of isolated components in the graphs of
these realizations.
In this section, we derive bounds on the number of diagonal blocks the
\(A\)-matrix of a system similar to a given system may have.
To do so, we primarily consider the class of systems with
block-diagonal \(A\)-matrices where each block is a companion matrix.
However, the upper bound we derive applies to all block-diagonal
realizations similar to the given system, no matter how their blocks
are structured.

\begin{definition}
A system \stdnamedsys is said to be a block-companion realization if
\(A\) is block-diagonal and each diagonal block \(A_i\) of \(A\) is a
companion matrix.
\end{definition}

In our first result, we state that there exists a lower bound on the number of
diagonal blocks of a block-companion realization similar to a given
system.
As we state in our second result, an upper bound on this number of
blocks also exists.

\begin{theorem}
  \label{thm:diag-block-lb}
  Let \stdnamedsys be a given linear system.
  Furthermore, let \(\phi_i\), \(1 \leq i \leq m\) be the irreducible
  polynomials that divide \dcharmat{A} and \(k_i\) the number of
  elementary divisors of \(A\) of the form \(\phi_i^l\).
  Then, every block-companion realization \(S'\) similar to \(S\) has
  at least \(k = \max_i k_i\) diagonal blocks.
\end{theorem}

\begin{theorem}
  \label{thm:diag-block-ub}
  Let \stdnamedsys be a given linear system and let \(d\) be its
  number of elementary divisors.
  Then, any block-companion realization \(S'\) similar to \(S\) has at
  most \(d\) diagonal blocks.
\end{theorem}

In the result below, we state that for any integer \(l\) between the
bounds we have previously indicated, we can find a block-companion
realization similar to \(S\) with \(l\) diagonal blocks.

\begin{theorem}
  \label{thm:diag-block-existence}
  Let \stdnamedsys be a given linear system and let \(k\) and \(d\) be
  the lower and upper bounds of Theorems \ref{thm:diag-block-lb} and
  \ref{thm:diag-block-ub} respectively.
  Then, for any integer \(l\) in the interval \([k,d]\), there exists
  a block-companion realization similar to \(S\) with \(l\) diagonal blocks.
\end{theorem}

As noted above, our main result is stated in terms of block-companion
realizations.
If we consider block-diagonal realizations whose blocks have arbitrary
shapes, we find that the upper bound given above still applies.
This fact is stated as a remark below.

\begin{remark}
  \label{rem:diag-block-ub}
  The upper bound of Theorem \ref{thm:diag-block-ub} also applies
  to arbitrary block-diagonal realizations, regardless of the shape of
  the blocks.
\end{remark}

While the upper bound of Theorem \ref{thm:diag-block-ub} remains
valid for arbitrary block-diagonal realizations, the lower bound does
not.
A counterexample is given below.
This example also indicates that this lower bound does not apply even
if we consider only minimal systems.

\begin{example}
 First, we will construct a system whose elementary divisors are
 \(\lambda -1\), with multiplicity 3, and \(\lambda-2\), with
 multiplicity 1.
 To do so, we note that the matrix \(A\) given below is in first
 natural normal form and has invariant polynomials \(i_1 =
 (\lambda-1)(\lambda-2)\), \(i_2 = i_3 = \lambda-1\).
 Let \namedsys{S}{A}{I_4}{I_4}{0}.
 Since \(B = C = I_4\), \(S\) is both controllable and observable.
 Thus, \(S\) is minimal, as required.
 Next, we consider \transys{S}{T}, with \(T\) given below.
 As is readily verified, \(\transys{S}{T} =
 \left(A_1,T,T^{-1},0\right)\), where \(A_1\) is given below.
 We notice that \(A_1\) has two diagonal blocks, even though the lower
 bound of Theorem \ref{thm:diag-block-lb} is three for \(S\).
  \begin{align*}
   A & = 
   \begin{bmatrix}
     0 & -2 & 0 & 0 \\
     1 & 3 & 0 & 0 \\
     0 & 0 & 1 & 0 \\
     0 & 0 & 0 & 1
   \end{bmatrix}\\
   T & = 
   \begin{bmatrix}
     -1 & 0 & 1 & 0 \\
     1 & 0 & 1 & 0 \\
     0 & 2 & 0 & 0 \\
     0 & 0 & 0 & 1
   \end{bmatrix}\\
   A_1 & = 
   \begin{bmatrix}
     \frac{1}{2} & \frac{1}{2} & 1 & 0 \\
     \frac{1}{2} & \frac{1}{2} & -1 & 0 \\
     -1 & 1 & 3 & 0 \\
     0 & 0 & 0 & 1
   \end{bmatrix}
 \end{align*}
\end{example}
\subsection{Structured systems and their graphs}
\label{subsec:main-structured-generic}
In many applications, we can assume that some variables are not
influenced by certain other variables.
In other words, we can assume that a number of edges
\((v_{i,1},v_{i,2})\),\(1 \leq i \leq l\), do not occur in the graph
of the system generating the input/output data under consideration.
This assumption implies that the graph \(G(S)\) of this system is a
subgraph of the graph \(G\) containing all the edges we have not
excluded, i.e. the graph \(G = (V,\{(v_1,v_2) \in V \times V | \neg
\exists i : v_1 = v_{i,1} \wedge v_2 = v_{i,2}\})\), where \(V\) is the
set of vertices of the graph \(G(S)\).
As we will see, many of the properties of systems in this class hold
for almost all systems in the class.
Furthermore, whether these properties hold for almost all systems can
be determined by examining the graph \(G\).
We call the class of systems \(S\) such that their graphs
\(G(S)\) are subgraphs of \(G\) a structured linear system.
Since the edges not present in \(G\) are absent in \(G(S)\) if and
only if the corresponding entries in the system matrices of \(S\) are
zero, a system \(S\) is a member of the structured linear system if
and only if those entries are zero.
Since the edges present in \(G\) may or may not be present in
\(G(S)\), the other entries of the system matrices of \(S\) are
unconstrained, i.e. they are free parameters.
We state a formal definition of a structured linear system below.

\begin{definition}
  \label{def:structured-matrix}
  Let \(T_{m \times n} = \{(i,j) | 1 \leq i \leq m \wedge 1 \leq j
  \leq n\}\) and \(T(M) \subset T_{m \times n}\).
  The structured matrix \(M\) corresponding to the subset \(T(M)\) is
  the subset of \(\mathbb{R}^{m \times n}\) given below.
  \[ M = \{M' \in \mathbb{R}^{m \times n} | \forall (i,j) \in T(M),
  M'_{ij} = 0\}\]
  The set of all \(m \times n\) structured matrices is denoted by
  \structuredspace{m}{n}.
\end{definition}

\begin{definition}
  \label{def:structured-ls}
  A structured linear system \(\mathbf{S}\) with \(n_u\) inputs, \(n_x\) state
  variables and \(n_y\) outputs is defined as a 4-tuple \((A,B,C,D)^S\),
  where \(A \in \structuredspace{n_x}{n_x}\), \(B
  \in \structuredspace{n_x}{n_u}\), \(C
  \in \structuredspace{n_y}{n_x}\) and \(D
  \in \structuredspace{n_y}{n_u}\).
  We say that \(\namedsys{S'}{A'}{B'}{C'}{D'} \in \mathbf{S}\) if \(A' \in A\),
  \(B' \in B\), \(C' \in C\) and \(D' \in D\).
\end{definition}

Above, we stated that the entries of a structured linear system that
were not set to zero are free parameters.
To make this precise, we introduce the standard parameterization of a
structured linear system.
In this parameterization, each entry in the system matrices that is
not set to zero corresponds to a single variable.

\begin{definition}
Let \(M \in \structuredspace{m}{n}\) and let \(T(M)\) and \(T_{m
  \times n}\) be given as in Definition \ref{def:structured-matrix}.
Furthermore, let \(\overline{T}(M) = T_{m \times n} \backslash T(M)\).
Order the elements of \(\overline{T}(M)\) according to the
lexicographic ordering 
\((i,j) \leq (i',j') \Leftrightarrow i < i' \vee (i = i' \wedge j \leq
j')\).
Let \(t_1, t_2, \cdots t_k\) be the elements of \(\overline{T}(M)\),
enumerated in the above order.
Then, the standard parameterization of the structured matrix \(M\) is
the function \(f_M : \mathbb{R}^k \rightarrow M\) given by \(f_M(v) = \left(e_{ij}(v)\right)\), with \(e_{ij}\) given
below.
\[e_{ij}(v) = 
\begin{cases}
  0 & \text{ if } (i,j) \in T(M) \\
  v_l & \text{ if } (i,j) = t_l
\end{cases}\]
The standard parameterization \(f_M\) is a bijection between
\(\mathbb{R}^k\) and \(M\).
\end{definition}

\begin{definition}
  Let \stdnamedsysstruct and let \(f_A : \mathbb{R}^{k_A} \rightarrow A\) be the standard
  parameterization of \(A\), with \(f_B\), \(k_B\), \(f_C\), \(k_C\),
  \(f_D\) and \(k_D\) defined similarly.
  Then, the standard parameterization of \(\mathbf{S}\) is the
  function \(f_\mathbf{S} : \mathbb{R}^{k_A+k_B+k_C+k_D} \rightarrow
  \mathbf{S}\) defined below.

  \[ f_\mathbf{S}\left(
  \begin{bmatrix}
    x_A \\
    x_B \\
    x_C \\
    x_D
  \end{bmatrix}\right) =
  \left(f_A(x_A),f_B(x_B),f_C(x_C),f_D(x_D)\right)\]

  The standard parameterization \(f_\mathbf{S}\) is a bijection
  between \(\mathbb{R}^{k_A+k_B+k_C+k_D}\) and \(\mathbf{S}\).
  We will denote \(f_\mathbf{S}(p)\) by \(\mathbf{S}_p\).
\end{definition}

A useful property of structured linear systems is that many of their
properties hold for almost all of systems within the class given by a
certain structured linear system.
Properties that hold for ``almost all'' systems are said to be
generic.
Formally, we say that a property \(P\) is generic if all parameter vectors
for which this property does not hold lie in the zero set of some
polynomial.
Equivalently, the set of parameter vectors for which \(P\) does not
hold is a subset of a proper variety as defined below.

\begin{definition}
  Let \(f\) be a polynomial function in \(n\) indeterminates with real
  coefficients.
  Then, the subset \(V_f = \{ p \in \mathbb{R}^n | f(p) = 0 \}\) is
  called the variety determined by \(f\).
  If \(V_f \not = \mathbb{R}^n\), \(V_f\) is said to be a proper variety.
\end{definition}
\begin{remark}
  In the remainder of this paper, we will denote the variety
  determined by \(f\) by \(V_f\). 
  Conversely, if we denote a variety by \(V_f\), the polynomial
  determining this variety is denoted by \(f\).
\end{remark}

\begin{definition}
  Let \(f_\mathbf{S} : \mathbb{R}^k \rightarrow \mathbf{S}\) be
  the standard parameterization of the structured linear system \(\mathbf{S}\).
  A property \(P\) is said to be generic for \(\mathbf{S}\) if for all \(p \in
  \mathbb{R}^k\) outside of a proper variety \(V\), \(P\) holds for \(\mathbf{S}_p\). 
\end{definition}

In our first result, we consider whether we can use input/output data
to uniquely identify a particular member of a structured linear
system.
In other words, we consider whether knowing that certain edges do not
occur in the system's graph is sufficient to allow us to identify the
system.
We say that a system can be uniquely identified using input/output
data if there does not exist a system with exactly the same
input/output behavior.
This is stated formally below.

\begin{definition}
  \stdnamedsys and \namedsys{S'}{A'}{B'}{C'}{D'} are equivalent if,
  given zero initial conditions, the outputs \(y_{S,k}\) of \(S\) and
  \(y_{S',k}\) of \(S'\) are identical for all \(k \geq 0\) for all
  input sequences \(u_k\).
\end{definition}

\begin{definition}
  Let \(\mathbf{S}\) be a structured linear system.
  Then, \(\mathbf{S}\) is identifiable for \(p\) if no parameter vector \(q\)
  exists such that \(p \not = q\) and \(\mathbf{S}_p\) and \(\mathbf{S}_q\) are
  equivalent.
\end{definition}

In the result below, we state that unless the system matrix \(C\)
consists only of fixed zeroes, the structured linear system \(S\) is
not generically identifiable.

\begin{theorem}
  \label{thm:structured-identifiability}
  Let \stdnamedsysstruct be a structured linear system such that at
  least one entry of \(C\) is not a fixed zero.
  Then, \(\mathbf{S}\) is not generically identifiable.
\end{theorem}

According to the result above, if a structured system
\stdnamedsysstruct is generically identifiable, all entries of the
matrix \(C\) are fixed zeroes.
In other words, for any parameter value \(p\), the output equation of
the system \(\mathbf{S}_p\) reads \(y_k = Du_k\).
Therefore, a structured system that is generically identifiable can
only represent very restricted input/output mappings.

Since all generically identifiable systems must have very limited
modeling power, all structured systems that will occur in practical
applications will not be generically identifiable.
This implies that in real applications, the knowledge we have of the
likely shape of the system's graph is insufficient to determine the
system parameters exactly.
In some cases, this knowledge may even be insufficient to uniquely
determine the system's structure.
To guarantee that the system's parameters can be found from
input/output data, additional knowledge about the system is necessary.

In the parameterizations we discussed in Section
\ref{subsec:main-cg-isomorphism}, we assumed that the system had a
particular shape and that it was minimal.
Together, these assumptions allowed us to obtain characterizations of
(CG-)isomorphism between systems.
In our second application of structured linear systems, we show how
the latter assumption, that is, that a system is minimal, implies that
a system's graph must satisfy a number of conditions.
This result can be used with any particular parameterization to study
the shape of the graph of a minimal system in this parameterization.
In addition, this result shows that the assumptions we make about the
properties a system has may affect its graph structure.

We begin by noting the following result, which states that a minimal
system must be a member of some generically minimal structured linear
system.

\begin{theorem}
  \label{thm:generic-minimal-is-existential}
  A structured linear system \(\mathbf{S}\) is generically minimal if and only
  if there exists a minimal system \(\mathbf{S}_p \in \mathbf{S}\).
\end{theorem}

This result implies that if a linear system \(S\) is minimal in the
usual sense, the structured linear system corresponding to its graph must
be generically minimal.
On the other hand, if the structured system linear system
corresponding to the graph of \(S\) is generically minimal, this does
not guarantee that \(S\) is minimal.

An ordinary linear system is minimal if and only if it is both
controllable and observable.
As we state below, the same holds generically for a structured linear
system.

\begin{theorem}
  \label{thm:generic-minimality}
  A structured linear system \(\mathbf{S}\) is generically minimal if and only
  if it is both generically controllable and generically observable.
\end{theorem}

As given by Dion et al.\cite{dion-et-al}, graph-theoretical conditions
for the generic controllability of a structured linear system exist.
Applying these conditions, we find the graph-theoretical
characterization of generic minimality stated below in Theorem \ref{thm:generic-minimality-graph}.
The graph of a structured linear system, to which this
characterization applies, is defined formally below.
Intuitively, an edge exists in this graph if and only if the
corresponding entry in the system matrices is not a fixed zero.

\begin{definition}
  The graph of a structured linear system \(\mathbf{S}\) is the graph of the system
  \(\mathbf{S}_v\), where \(v\) is the vector such that \(v_i = 1\) for all \(i\).
\end{definition}

Below, we define the concepts used in the statement of Theorem
\ref{thm:generic-minimality-graph}.
Similar definitions of these concepts are given by Dion et
al.\cite{dion-et-al}.

\begin{definition}
  Let \stdnamedsysstruct be a structured linear system.
  A finite sequence of vertices \(v_1 v_2\cdots v_k\) in \(G(\mathbf{S})\) such
  that no two vertices \(v_i\) and \(v_j\) are equal except possibly \(v_1\) and
  \(v_k\) is called a path of length \(k\).
  If \(v_1 = v_k\), the path is called a cycle. Otherwise, the path is
  a simple path.
  Furthermore, if \(v_1 = u_i(G(\mathbf{S}))\) for some \(i\), the path is
  called \(U\)-rooted.
  Similarly, if \(v_k = y_j(G(\mathbf{S}))\) for some \(j\), the path is said
  to be \(Y\)-topped.
\end{definition}

\begin{definition}
  A path \(u_1 u_2 \cdots u_k\) is said to cover the vertices \(u_1\)
  through \(u_k\).
  Two paths \(p = u_1 u_2 \cdots u_k\) and \(q = u'_1 u'_2 \cdots
  u'_l\) are said to be disjoint if no vertex \(w\) exists such that
  both \(p\) and \(q\) cover \(w\).

  A set of simple paths such that any two of them are disjoint is
  called a path family.
  If every element of the family is \(U\)-rooted, the family is said
  to be \(U\)-rooted.
  Similarly, a path family of which every member is \(Y\)-topped is
  called \(Y\)-topped.  
  A set of cycles such that any two of them are disjoint is called a
  cycle family.

  The union of two path or cycle families \(F_1\) and \(F_2\) is said
  to be disjoint if for all \(p_1 \in F_1\) and \(p_2 \in F_2\),
  \(p_1\) and \(p_2\) are disjoint.
\end{definition}

\begin{theorem}
  \label{thm:generic-minimality-graph}
  A structured linear system \(\mathbf{S}\) is generically minimal if and only
  if its graph \(G(\mathbf{S})\) satisfies the following conditions.
  \begin{enumerate}
  \item Every vertex \(x_i(G(\mathbf{S}))\) is the end vertex of a \(U\)-rooted
    path in \(G(\mathbf{S})\)
  \item There exists a disjoint union of a \(U\)-rooted path family
    and a cycle family in \(G(\mathbf{S})\) that covers all the vertices
    \(x_i(G(\mathbf{S}))\)
  \item Every vertex \(x_i(G(\mathbf{S}))\) is the first vertex of a
    \(Y\)-topped path in \(G(\mathbf{S})\)
  \item There exists a disjoint union of a \(Y\)-topped path family
    and a cycle family in \(G(\mathbf{S})\) that covers all the vertices \(x_i(G(\mathbf{S}))\).
  \end{enumerate}
\end{theorem}

As we have said previously, we can use the graph-theoretical
characterization of generic minimality given by the previous result to
determine if a given linear system may be minimal.
This result is stated formally below.

\begin{theorem}
  \label{thm:minimality-graph-necessary}
  If a linear system \(S\) is minimal, its graph must satisfy the
  conditions of Theorem \ref{thm:generic-minimality-graph}.
\end{theorem}

\section{Proofs of the main results}
\label{sec:proofs-main-results}
\subsection{Graph isomorphism and its inadequacy}
We begin by stating the proof of Theorem
\ref{thm:diagonalizable-non-isomorphism}.

\begin{proof}[Proof of Theorem \ref{thm:diagonalizable-non-isomorphism}]
Since \(A\) is diagonalizable, there exists an invertible matrix \(T\)
such that \(A = TMT^{-1}\), where \(M\) is a diagonal matrix.
Equivalently, \(M = T^{-1}AT\).
This implies that \(\transys{S}{T^{-1}} = (M,B',C',D')\), for some
\(B'\), \(C'\) and \(D'\).
Since \(M\) is diagonal, no edge of the form
\((x_i(G(\transys{S}{T^{-1}})),x_j(G(\transys{S}{T^{-1}})))\) exists in
\(G(\transys{S}{T^{-1}})\), for distinct \(i\) and \(j\).
However, since \(A\) was not diagonal, at least one edge of this form
exists in \(G(S)\).
An isomorphism must preserve this edge, which is impossible.
Therefore, \(G(S)\) and \(G(\transys{S}{T^{-1}})\) are non-isomorphic,
as claimed.
\end{proof}

In order to prove Theorem \ref{thm:gi-n-subgroup}, we first take a
closer look at the type-restricted isomorphism relation.
We formally re-define this relation below.

\begin{definition}
\label{def:tr-isomorphism}
The systems \(S\) and \(S'\) are said to be type-restricted
isomorphic, denoted \(S \simeq S'\), if a type-restricted isomorphism
exists between \(G(S)\) and \(G(S')\).
\end{definition}

As we state below, the relation defined above is an equivalence
relation.
The proof of this lemma is deferred to an appendix.

\begin{lemma}
  \label{lem:tr-isomorphism-equiv}
  Type-restricted isomorphism (written \(S \simeq S'\)), is an
  equivalence relation, that is:
  \begin{enumerate}
  \item \(S \simeq S\), for all \(S\)
  \item If \(S \simeq S'\), then \(S' \simeq S\)
  \item If \(S \simeq S'\) and \(S' \simeq S''\), then \(S \simeq S''\)
  \end{enumerate}
\end{lemma}

Using Lemma \ref{lem:tr-isomorphism-equiv}, we can now state the proof
of Theorem \ref{thm:gi-n-subgroup}.

\begin{proof}[Proof of Theorem \ref{thm:gi-n-subgroup}]
To prove the first condition, notice that \(\transys{S}{I_n} = S\) for
all \(S\).
Therefore, \(\transys{S}{I_n} \simeq S\).

Next, let \(T,T' \in GI(n)\) and let \(S\) be some realization.
Then, \(\transys{S}{T} \simeq S\), by the definition of \(GI(n)\).
For the same reason, \(\transys{\transys{S}{T}}{T'} \simeq
\transys{S}{T}\), and so \(\transys{\transys{S}{T}}{T'} \simeq S\).
Noticing that \(\transys{\transys{S}{T}}{T'} = \transys{S}{T'T}\), we
find that \(T'T \in GI(n)\), as required.

Finally, let \(T \in GI(n)\).
Then, let \(S' = \transys{S}{T^{-1}}\), for some \(S\).
It follows that \(S = \transys{S'}{T}\).
Therefore, \(S \simeq S'\), by the definition of \(GI(n)\).
Thus, \(S' \simeq S\) and so \(T^{-1} \in GI(n)\).
\end{proof}

Next, we state the proof of Theorem \ref{thm:diag-gi-n}.

\begin{proof}[Proof of Theorem \ref{thm:diag-gi-n}]
  Let \(M\) be given and \(\stdnamedsys\) be an arbitrary
  realization.
  Then, \(\transys{S}{M} = (MAM^{-1},MB,CM^{-1},D)\).
  Since \(M\) is diagonal, so is \(M^{-1}\), and \(M^{-1} = \{m_i^{-1}\}\).
  Notice that if \(L = (L_{ij})\),\footnote{\((L_{ij})\) denotes the
  matrix whose elements are \(L_{ij}\)} \(ML = (m_iL_{ij})\).
  Similarly, \(LM^{-1} = (L_{ij}m_j^{-1})\).
  Therefore, \(\transys{S}{M} =
  ((m_iA_{ij}m_j^{-1}),(m_iB_{ij}),(C_{ij}m_j^{-1}),D)\).
  Since both \(m_i\) and \(m_i^{-1}\) are clearly non-zero for all
  \(i\), each entry of each matrix in \transys{S}{M} is non-zero if
  and only if the corresponding entry in \(S\) is non-zero.
  Therefore, the graph of \(S\) is the same as that of \transys{S}{M},
  and so \(S \simeq \transys{S}{M}\).
  Since \(S\) was arbitrary, \(M \in GI(n)\), as claimed.
\end{proof}

To prove that permutation matrices are elements of \(GI(n)\), we
recall a number of properties of such matrices below.
The first of these properties is proved in our appendix.
The second follows straightforwardly from matrix multiplication.

\begin{lemma}
\label{lem:permutation-orthogonal}
A permutation matrix \(P(e_1,e_2,\cdots,e_n)\) is orthogonal.
\end{lemma}

\begin{lemma}
\label{lem:permutation-products}
For a matrix \(A = (A_{ij})\) and a permutation matrix
\(P(e_1,e_2,\cdots,e_n)\), \(P(e_1,e_2,\cdots,e_n)A = (A_{e_ij})\) and
\(AP(e_1,e_2,\cdots,e_n)^T = (A_{ie_j})\), whenever these products exist.
\end{lemma}

Below, we show that permutation matrices are indeed elements of
\(GI(n)\), proving Theorem \ref{thm:permutation-gi-n}.

\begin{proof}[Proof of Theorem \ref{thm:permutation-gi-n}]
Let \namedsys{S}{(A_{ij})}{(B_{ij})}{(C_{ij})}{(D_{ij})} and let \(P =
P(e_1,e_2,\cdots,e_n)\) be an arbitrary permutation matrix.
By Lemma \ref{lem:permutation-products}, \(P(B_{ij}) = (B_{e_ij})\) and \((C_{ij})P^T =
(C_{ie_j})\).
Similarly, \(P(A_{ij}) = (A_{e_ij})\), and so \(P(A_{ij})P^T =
(A_{e_ie_j})\).
Thus, \(\transys{S}{P} =
\left((A_{e_ie_j}),(B_{e_ij}),(C_{ie_j}),D\right)\).
We claim that \(\phi\) defined below is a type-restricted isomorphism
between \(S\) and \(\transys{S}{P}\). 
\[ \phi(v) = \begin{cases}
u_i(G(\transys{S}{P})) & \text{if } v = u_i(G(S)) \\
y_i(G(\transys{S}{P})) & \text{if } v = y_i(G(S)) \\
x_i(G(\transys{S}{P})) & \text{if } v = x_{e_i}(G(S))
\end{cases}
\]
From the definition of \(\phi\), it is clear that \(\phi\) is
type-restricted.
Since the integers \(e_1,e_2,\cdots,e_n\) are a permutation of the
integers from 1 through \(n\), it is also clear that \(\phi\) is a
bijection.
Thus, we only need to verify that \(\phi\) preserves edges.
To do so, we consider the four types of edges given in Corollary
\ref{cor:associated-graph}.

First, we have edges of the form \((x_{e_j},x_{e_i})\).
Then, \(\phi(x_{e_j}) = x_j\) and \(\phi(x_{e_i}) = x_i\).
The edge \((x_j,x_i)\) exists in \(G(\transys{S}{P})\) if and
only if \((PAP^T)_{ij} = (A_{e_ie_j})_{ij} = 
A_{e_ie_j} \not = 0\).
Thus, this edge exists if and only if \((x_{e_j},x_{e_i})\) is an edge of
\(G(S)\).

Second, consider edges of the form \((u_j,x_{e_i})\), which exist in
\(G(S)\) if and only if \(B_{e_ij} \not = 0\).
By definition, \((\phi(u_j),\phi(x_{e_i})) = (u_j,x_i)\).
This edge exists in \(G(\transys{S}{P})\) if and only if \((PB)_{ij}
= (B_{e_ij})_{ij} = B_{e_ij} \not = 0\).
Thus, this edge is conserved.

Third, consider an edge of the form \((x_{e_j},y_i)\), which is an edge of
\(G(S)\) if and only if \(C_{ie_j} \not =  0\).
Using \(\phi\) again, \((\phi(x_{e_j}),\phi(y_i)) = (x_j,y_i)\).
This edge exists in \(G(\transys{S}{P})\) if and only if
\((CP^T)_{ij} = (C_{ie_j})_{ij} = C_{ie_j} \not = 0\).
Therefore, edges of this form are also conserved.

Finally, all edges of the form \((u_i,y_j)\) correspond to edges
\((u_i,y_j)\), since outputs and inputs are not permuted by \(\phi\).
Since both \(S\) and \(\transys{S}{P}\) have the same matrix \(D\),
these edges are clearly conserved.
Thus, \(\phi\) is a type-restricted isomorphism.
Since \(S\) was arbitrary, \(P \in GI(n)\).
\end{proof}
\subsection{Condensed-graph isomorphism}
Below, we formally prove Theorem \ref{thm:condensed-diag-counterexample}.

\begin{proof}[Proof of Theorem \ref{thm:condensed-diag-counterexample}]
  First, we note that since \(D\) is diagonal, \(G(\transys{S}{T})\)
  has no components with more than one element.
  Thus, if \(G(S)\) has any components with more than one element, \(CG(S)\) and
  \(CG(\transys{S}{T})\) are of different order and so cannot be
  isomorphic.
  Second, if \(G(S)\) has no such components, \(CG(S)\) and
  \(CG(\transys{S}{T})\) are identical to \(G(S)\) and
  \(G(\transys{S}{T})\) respectively.
  We have shown previously that no type-restricted isomorphism between
  these graphs can exist.
\end{proof}

We will now prove Theorem \ref{thm:1-n-components} in two steps.
In the first step, we construct a transfer function such that the
realization in observable canonical form of this transfer function has 1
state variable component in its condensed graph.
Then, we show that we can construct a transfer function of this type
where the resulting realization has a diagonalizable \(A\)-matrix.

\begin{lemma}
Let \(H(s) = \frac{b_0\lambda^n
  +b_1\lambda^{n-1}+\cdots+b_n} {\lambda^n +
  a_1\lambda^{n-1}+\cdots+a_n}\) be a given transfer function.
Then, the condensed-graph of the observable canonical form of \(H(s)\)
contains exactly one state-variable component if \(a_n \not = 0\).
\end{lemma}
\begin{proof}
\label{ex:observable-cf}
Consider a transfer function \(H(s) = \frac{b_0\lambda^n
  +b_1\lambda^{n-1}+\cdots+b_n} {\lambda^n +
  a_1\lambda^{n-1}+\cdots+a_n}\).
The observable canonical form of \(H(s)\), as described by Vaccaro
\cite[Ch. 3.3.2]{vaccaro}, is given by \stdnamedsys, where \(A\), \(B\),\(C\) and
\(D\) are given below.
We notice that for all values of the coefficients \(a_i\), the edges
\((x_i,x_{i-1})\) for \(2 \leq i \leq n\) exist in \(G(S)\).
If \(a_n \not = 0\), the edge \((x_1,x_n)\) will also exist in
\(G(S)\), completing a cycle containing all state variables.
Thus, as claimed, \(G(S)\) will then have a single component
containing all state variables.

\begin{align*}
A & = \begin{bmatrix}
-a_1 & 1 & 0 & 0 & \cdots & 0 \\
-a_2 & 0 & 1 & 0 & \cdots & 0 \\
-a_3 & 0 & 0 & 1 & \cdots & 0 \\
\vdots & \vdots & \vdots & \vdots & \ddots & \vdots \\
-a_{n-1} & 0 & 0 & 0 & \cdots & 1 \\
-a_n & 0 & 0 & 0 & \cdots & 0
\end{bmatrix}
& B & = 
\begin{bmatrix}
  b_1-a_1b_0 \\
  b_2-a_2b_0 \\
  b_3-a_3b_0 \\
  \vdots \\
  b_n -a_nb_0
\end{bmatrix} \\ C & = 
\begin{bmatrix}
  1 & 0 & \cdots & 0
\end{bmatrix} 
& D  & = 
\begin{bmatrix}
  b_0
\end{bmatrix}
\end{align*}
\end{proof}

To complete the proof of Theorem \ref{thm:1-n-components}, we
construct a transfer function satisfying the conditions of the
previous lemma such that the observable canonical form has a
diagonalizable \(A\)-matrix.
This construction is given below.

\begin{proof}[Proof of Theorem \ref{thm:1-n-components}]
  Let \(H(s) = \frac{1}{\prod_{i=1}^n(s-i)}\).
  Since \((0-i)\) is non-zero for any \(1 \leq i \leq n\), the
  coefficient \(a_n\) of the denominator of \(H(s)\) is non-zero.
  Let \(S\) be the observable canonical form of \(H(s)\).
  As we have shown previously, \(G(S)\) contains a component
  containing all state variables.
  
  Let \(A\) be the \(A\)-matrix of \(S\).
  Since the eigenvalues of \(A\) are the poles of \(H(s)\), these
  eigenvalues are the integers from 1 through \(n\).
  As \(A\) has \(n\) distinct eigenvalues, \(A\) is diagonalizable.
  Thus, we can find a matrix \(T\) such that \(G(\transys{S}{T})\) has \(n\)
  components containing state variables.
\end{proof}

In order to prove our isomorphism results, we introduce a number of concepts
that we have found useful.
Both of these concepts are sets of state variables.
We recall from Definition \ref{def:associated-graph} that the state
variables in the graph \(G(S)\) of some system \(S\) are the vertices \(x_i(G(S))\).
The first concept is that of a trap, i.e.\ a set of state variables with no
edges to any vertex not in the set itself.
The second is the similar concept of an unreachable set, i.e.\ a set of
state variables such that no vertex outside of the set itself has an
edge to any of the variables in the set.
These concepts are defined formally below.

\begin{definition}
\label{def:trap}
A non-empty set of state variables \(S\) in the graph \(G\) of some system is called a
trap if there does not exist a vertex \( x \in S \) and a vertex \( y
\not \in S \) such that \((x,y)\) is an edge of G.
\end{definition}

\begin{definition}
\label{def:unreachable-set}
A non-empty set of state variables \(S\) in the graph \(G\) of some system is
called an unreachable set if there does not exist a vertex \(x \in S\)
and a vertex \(y \not \in S \) such that \((y,x)\) is an edge of G.
\end{definition}

A useful property of traps and unreachable sets is that any system
whose graph contains either an unreachable set or a trap is
non-minimal.
We state this below in two lemmas.
The proof of the first of these is in an appendix, the second has a
proof very similar to that of the first.

\begin{lemma}
\label{lem:trap-non-observability}
A realization \stdnamedsys whose graph contains a trap is not
observable, that is, its observability matrix \(\mathfrak{O} = \left[
\begin{array}{c}
C \\
CA \\
\vdots \\
CA^{n-1}
\end{array}\right]\) has rank less than n.
\end{lemma}

\begin{lemma}
\label{lem:unreachable-set-uncontrollability}
A realization \((A,B,C,D)\) whose graph contains an unreachable set is
uncontrollable, that is, the rank of its controllability matrix
\(\mathfrak{C} = \left[\begin{array}{cccc}
B & AB & \cdots & A^{n-1}B
\end{array}\right]\) is less than \(n\).
\end{lemma}

\begin{remark}
  A result similar to Lemma \ref{lem:unreachable-set-uncontrollability}
  is described by Dion et al.\cite{dion-et-al}.
  Though we will not formally prove this here, the presence of an
  unreachable set in a system's graph is equivalent to the system's
  being in Form I as given by Dion et al.\cite{dion-et-al}.
\end{remark}

We now note that if, in \stdnamedsys, the matrix \(A\) is diagonal, no
edges between distinct state variables exist.
Thus, each state variable is a trap if it has no edge to the output
\(y\), and an unreachable set if it has no edge from the input.
We state this formally below.

\begin{observation}
In the graph of a minimal diagonal realization
, each state variable has an edge from the input
and an edge to the output.
\end{observation}

We now note that an isomorphism between diagonal realizations will
preserve edges involving a state variable \(x_i\) if and only if it
preserves edges of the form \((x_i,x_i)\).
This is stated formally below.

\begin{observation}
\label{obs:diag-isomorphism-sl}
An isomorphism \(\phi\) between the graphs \(G = (V,E)\) and \(G' = (V',E')\) of two
diagonal realizations may have \(\phi(x_i) = x_j\) if and only if either
\((x_i,x_i) \in E \wedge (x_j,x_j) \in E'\) or \((x_i,x_i) \not \in E
\wedge (x_j,x_j) \not \in E'\)
\end{observation}

Using the above observations, we will now prove Theorem \ref{thm:diag-isomorphism}.

\begin{proof}[Proof of Theorem \ref{thm:diag-isomorphism}]
To prove necessity, let \(G(S)\) and \(G(S')\) be isomorphic.
Then, since \((u,y)\) is an edge of either both \(G(S)\) and
\(G(S')\) or of neither of these graphs, either both \(D\) and \(D'\)
are non-zero or \(D = D' = 0\).
Since any isomorphism between \(G(S)\) and \(G(S')\) must conserve
edges of the form \((x_i,x_i)\), \(G(S)\) and \(G(S')\) have equal
numbers of vertices for which such an edge exists.
Since each such vertex corresponds to a non-zero diagonal entry of the
matrix \(A\) or \(A'\), the numbers of such entries are equal.

To prove sufficiency, we construct a bijection \(\phi\).
Let \(x_{11},x_{12},\cdots,x_{1m_1}\) be the state variable vertices in
\(G(S)\) for which an edge of the form \((x_i,x_i)\) exists, that is,  the
state variable vertices corresponding to non-zero diagonal entries of \(A\).
Let \(x_{21},x_{22},\cdots,x_{2m_2}\) be the other state variable
vertices of \(G(S)\).
Define \(x'_{11},x'_{12},\cdots,x'_{1m_1}\) and
\(x'_{21},x'_{22},\cdots,x'_{1m_2}\) similarly for \(G(S')\).
Then, we claim that the function \(\phi\) defined below is an
isomorphism between \(G(S)\) and \(G(S')\).

\[\phi(v) = \begin{cases}
u_1(G(S')) &\text{ if } v = u_1(G(S)) \\
y_1(G(S')) &\text{ if } v = y_1(G(S)) \\
x'_{1i} & \text{ if } v = x_{1i} \\
x'_{2i} & \text{ if } v = x_{2i}
\end{cases}\]

\(\phi\) is clearly a bijection.
It is also clear from the definition of \(\phi\) that \(\phi\) is
type-restricted.
Thus, we only need to verify that \(\phi\) preserves edges.
Since \(S\) is minimal, we know that for any state
variable \(x_i\), \((u,x_i)\) and \((x_i,y)\) are edges of \(G(S)\).
For the same reason, we know the same holds for \(\phi(x_i)\).
Since both \(A\) and \(A'\) are diagonal, no edges between distinct
state variables exist in either \(G(S)\) or \(G(S')\).
By the definition of \(\phi\), it is clear that \(\phi\) also
conserves edges of the form \((x_i,x_i)\).
Finally, since \(D\) and \(D'\) are either both zero or both non-zero,
\((u,y)\) is either an edge of both \(G(S)\) and \(G(S')\) or of
neither of these graphs.
\end{proof}

To prove Theorem \ref{thm:first-natural-companion}, we will use the
following lemma.
The proof of this result is deferred to an appendix.

\begin{lemma}
\label{lem:minpoly-non-minimality}
If, for a SISO realization \((A,B,C,D)\), the minimal polynomial \(i_1\) of
\(A\) is not equal to \(\dcharmat{A}\) the realization is non-minimal.
\end{lemma}

Using the above lemma, we state the proof of Theorem
\ref{thm:first-natural-companion} below.

\begin{proof}[Proof of Theorem \ref{thm:first-natural-companion}]
Using the above lemma, we see that if \stdnamedsys is a minimal SISO
realization , the minimal polynomial \(i_1\) of \(A\) coincides with
\dcharmat{A}.
Since the product of all minimal polynomials \(i_j\) of \(A\) also
coincides with \dcharmat{A}, \(i_j = 1\) for \(j > 1\).
Thus, the first natural normal form of \(A\) consists of a single
diagonal block that is the companion matrix of \(i_1 = \dcharmat{A}\).
Therefore, we have established Theorem \ref{thm:first-natural-companion}.  
\end{proof}

To prove our isomorphism result for realizations with \(A\)-matrices
in second natural normal form, we begin by considering the graph of
each block of the second natural normal form.
As stated below, these graphs either have Hamiltonian cycles or are
directed paths.

\begin{lemma}
\label{lem:elem-div-structure}
For an irreducible polynomial \(\phi(\lambda)\), the graph of the companion matrix
for polynomials of the form \(\phi(\lambda)^k\), where \(k > 0\),
satisfies one of the following conditions:
\begin{enumerate}
\item \(\phi(0) \not = 0\) and the graph is Hamiltonian
\item \(\phi(0) = 0\) and the graph is a directed path.
\end{enumerate}
\end{lemma}

Since the graph of each block of \(A\) is a subgraph of the graph of
\stdnamedsys, we have the following corollary of Lemma \ref{lem:elem-div-structure}.

\begin{corollary}
If \(S = (A,B,C,D)\), where \(A\) is in second natural normal form,
each block of \(A\) corresponds to either a Hamiltonian component or a
directed path in \(G(S)\).
\end{corollary}

By applying the ideas of traps and unreachable sets to realizations in
second natural normal form, we find the following lemma.

\begin{lemma}
\label{lem:second-natural-component-conditions}
If \(S = (A,B,C,D)\) is a minimal SISO realization  and \(A\) is in
second natural normal form, each  component consisting of state
variables in \(G(S)\) satisfies one of the following two conditions.
Furthermore, these conditions are sufficient to guarantee that
\(G(S)\) contains neither traps nor unreachable sets.
\begin{enumerate}
\item The component is Hamiltonian and state variables \(x_i\) and
  \(x_j\) in the component exist such that \((u,x_i)\) and \((x_j,y)\)
  are both edges of \(G(S)\)
\item The component is a directed path consisting of state variables
  \(x_i,x_{i+1},\cdots,x_{i+m}\) and both \((u,x_i)\) and
  \((x_{i+m},y)\) are edges of \(G(S)\).
\end{enumerate}
\end{lemma}

Additionally, we obtain the below lemma and its corollary by applying Lemma
\ref{lem:minpoly-non-minimality} to realizations in second natural
normal form.

\begin{lemma}
\label{lem:minimal-irreducible-poly-uniqueness}
If \(S = (A,B,C,D)\) is a SISO realization, where \(A\) has two elementary divisors of the
form \(\phi(\lambda)^k\), for some irreducible \(\phi(\lambda)\),
\(S\) is non-minimal.
\end{lemma}

\begin{corollary}
If \(S = (A,B,C,D)\) is a minimal SISO realization and \(A\) is in second natural normal
form, the number of blocks in \(A\) equals the number of distinct
irreducible polynomials that divide \dcharmat{A}.
\end{corollary}

Using the above lemmas, we now state the proof of Theorem \ref{thm:second-natural-isomorphism}.

\begin{proof}[Proof of Theorem \ref{thm:second-natural-isomorphism}]
  To prove necessity, suppose \(G(S_1)\) and \(G(S_2)\) are
  CG-isomorphic.
  Then, these graphs must have the same number of state-variable
  components.
  Since 0 is not an eigenvalue of either \(A_1\) or \(A_2\), each
  block in \(A_1\) and \(A_2\) corresponds to a Hamiltonian component,
  by Lemma \ref{lem:elem-div-structure}.
  Thus, the number of such blocks must be equal.
  The first condition then follows from our corollary above.
  To show the second condition, notice that the edge \((u,y)\) in
  \(G(S_1)\) must exist if and only if the same edge exists in
  \(G(S_2)\).
  The second condition is obtained by restating this in terms
  of \(D_1\) and \(D_2\).

  To prove sufficiency, suppose the conditions are satisfied.
  Then, let \(c_1,c_2,\cdots,c_m\) and \(c'_1,c'_2,\cdots,c'_{m'}\) be
  the state-variable components of \(G(S_1)\) and \(G(S_2)\)
  respectively.
  By our corollary above, \(A_1\) and \(A_2\) have an equal number of
  diagonal blocks.
  As used above, this implies that \(m = m'\).
  We will now show that \(\phi\), as defined below, is a
  CG-isomorphism between \(G(S_1)\) and \(G(S_2)\).
  \[ \phi(c) = \begin{cases}
    u_1(CG(S_2)) & \text{if } c = u_1(CG(S_1)) \\
    y_1(CG(S_2)) & \text{if } c = y_1(CG(S_1)) \\
    c'_i & \text{if } c = c_i
    \end{cases}\]
  Clearly, \(\phi\) is type-restricted and, since \(m = m'\), \(\phi\)
  is a bijection.
  Thus, we verify that \(\phi\) conserves edges.
  Since each diagonal block in both \(A_1\) and \(A_2\) is a companion
  matrix corresponding to a polynomial with no root equal to zero,
  each component is either non-trivial and Hamiltonian or trivial and
  corresponding to a non-zero \((1 \times 1)\) matrix.
  Therefore, \((c_i,c_i)\) is an edge of \(CG(S_1)\) for all \(i\),
  and similarly for \((c'_i,c'_i)\) in \(CG(S_2)\).
  Since each component corresponds to a diagonal block, no edge of the
  form \((c_i,c_j)\) for distinct \(i\) and \(j\) exists in
  \(CG(S_1)\) and similarly in \(CG(S_2)\).
  Furthermore, since both realizations are minimal, \((u,c_i)\) and
  \((c_i,y)\) are edges of \(CG(S_1)\) and similarly in \(CG(S_2)\).
  Finally, \((u,y)\) is an edge of \(CG(S_1)\) if and only if \(D_1
  \not = 0\) and similarly for \(CG(S_2)\) and \(D_2\).
  Since \(D_1\) and \(D_2\) are either both zero or both non-zero,
  \(\phi\) preserves this edge.
\end{proof}
\subsection{Components of condensed graphs}

To prove our main results in this subsection, we will need the
following theorem from Gantmacher \cite[Ch.~6,Thm.~5]{gantmacher}.
This result allows us to find the elementary divisors of a block-diagonal
matrix using those of the diagonal blocks.

\begin{theorem}
\label{thm:elem-div-diagonal}
Let \(A = \{A_i\}\) \footnote{Recall that \(\{A_i\}\) denotes the
  block-diagonal matrix consisting of the blocks \(A_i\)}.
Then, \(A\) has all the elementary divisors of the \(A_i\), and no
others.
\end{theorem}

We begin by proving Theorem \ref{thm:diag-block-lb}, i.e. that there
exists a certain lower bound on the number of diagonal blocks in a
block-companion realization similar to a given system.
To do so, we recall from our definition of a companion matrix that
each companion matrix \(L\) has only one invariant polynomial not equal to
\(1\).
This implies that for any irreducible polynomial \(\phi\) that divides
\dcharmat{L}, \(L\) has exactly one elementary divisor of the form \(\phi^k\).

\begin{proof}[Proof of Theorem \ref{thm:diag-block-lb}]
Let \namedsys{S'}{A'}{B'}{C'}{D'} be an arbitrary block-diagonal
realization consisting of companion blocks that is similar to \(S\).
Each block of \(A'\) is a companion matrix and so cannot have two
elementary divisors corresponding to the same irreducible polynomial.
Since \(A'\) is similar to \(A\), \(A'\) has the same elementary
divisors as \(A\).
Thus, there exists an integer \(i\) such that \(k_i = k\) and \(A'\)
has \(k_i\) elementary divisors of the form \(\phi_i^l\).
Since each block of \(A'\) has at most 1 elementary divisor of this
form, \(A'\) must have at least \(k\) such blocks.
\end{proof}

In a block-companion realization, each block of the matrix \(A\)
corresponds to some nonzero number of elementary divisors of \(A\).
Thus, the maximum number of diagonal blocks is the number of
elementary divisors of \(A\).
This argument establishes the upper bound of Theorem
\ref{thm:diag-block-ub}.
In our discussion above, we did not use the assumption that the diagonal blocks of
\(A\) were companion matrices.
Thus, we have also proven Remark \ref{rem:diag-block-ub}.

The proof of the remaining result, Theorem
\ref{thm:diag-block-existence}, is split into two parts.
First, we show that the bounds given by the previous two results are sharp.
In other words, there exist realizations with a number of blocks equal
to the lower bound and realizations with a number of blocks equal to
the upper bound.
These results are stated below as lemmas.
The proofs of these lemmas are given in Appendix \ref{app:tech-condensed-components}.

\begin{lemma}
\label{lem:diag-block-lb-sharpness}
Let \stdnamedsys and \(k\) be the bound of Theorem
\ref{thm:diag-block-lb}.
Then, there exists a block-companion realization \(S'\) similar to \(S\) such that
\(S'\) has \(k\) diagonal blocks.
\end{lemma}

\begin{lemma}
\label{lem:diag-block-ub-sharpness}
Let \(A\) have \(n\) elementary divisors.
Then, there exists a block-companion realization similar to
\stdnamedsys with \(n\) diagonal blocks.
\end{lemma}

In the following lemma, we state that if we have a block-companion realization with
\(l\) blocks, we can rearrange the elementary divisors of these blocks
to find a block-companion realization with \(l+1\) blocks, provided
one exists.
This lemma allows us to prove the remainder of Theorem
\ref{thm:diag-block-existence}, namely that a block-companion realization
with \(l\) blocks exists for any integers \(l\) between the bounds indicated.

\begin{lemma}
\label{lem:diag-block-inc}
Let \stdnamedsys  be a block-companion realization such that \(A\) has
\(n\) elementary divisors and \(l\) diagonal blocks.
Then, either \(l = n\) or there exists a block-companion realization
similar to \(S\) with \(l+1\) diagonal blocks.
\end{lemma}

\begin{proof}[Proof of the remainder of Theorem \ref{thm:diag-block-existence}]
By Lemma \ref{lem:diag-block-lb-sharpness}, there exists a block-companion realization
similar to \(S\) with \(k\) diagonal blocks.
By applying Lemma \ref{lem:diag-block-inc} \(l - k\) times, we can
find a sequence of realizations \(S_1,S_2,\cdots,S_{l-k}\) such that
\(S_i\) is a block-companion realization similar to \(S\) and \(S_i\)
has \(k + i\) diagonal blocks.
Thus, \(S_{l-k}\) is the required realization.
\end{proof}

\subsection{Structured systems and their graphs}
The first result we will prove is Theorem \ref{thm:structured-identifiability}.
To do so, we begin by stating the fact that \transys{S}{M} and \(S\)
have the same input/output relation in terms of our definition of
equivalence.
This fact is well-known and is stated in many textbooks on linear
systems, such as the book by Vaccaro \cite[Section 3.3.4]{vaccaro}.
We reproduce a formal proof of this result in the appendix.

\begin{lemma}
  \label{lem:trans-equivalence}
  \stdnamedsys and \transys{S}{M} are equivalent for all \(M\).
\end{lemma}

The next step of our proof is to show that a structured linear system
satisfying the condition of Theorem
\ref{thm:structured-identifiability} is generically not identifiable.
We note that this result is subtly different from our theorem, which
requires us to establish that this kind of structured linear system is
not generically identifiable.
Still, the below lemma constitutes a major part of our proof.
Hence, its proof is stated here in full.

\begin{lemma}
  \label{lem:generic-non-identifiable}
  Let \stdnamedsysstruct be a structured linear system such that at
  least one entry of \(C\) is not a fixed zero.
  Then, \(\mathbf{S}\) is generically not identifiable.
\end{lemma}
\begin{proof}
  Let \(x_i\) be the parameter corresponding to the entry of \(C\)
  that is not a fixed zero.
  Let \(V_f\) be the variety corresponding to \(f(x) = x_i\).
  We claim that for \(p \not \in V_f\), \(\mathbf{S}\) is not identifiable.
  Thus, if \(\mathbf{S}\) is identifiable for \(q\), \(q \in V_f\), from which
  our result follows.
  Let \(p \not \in V_f\) and \(\mathbf{S}_p = (A_p,B_p,C_p,D_p)\).
  Then, \(\transys{\mathbf{S}_p}{2I} = (A_p,2B_p,\frac{1}{2}C_p,D_p)\).
  Thus, every entry of \transys{\mathbf{S}_p}{2I} corresponding to a fixed zero
  in \(\mathbf{S}\) is zero, since these entries are zero in \(\mathbf{S}_p\).
  Therefore, \(\transys{\mathbf{S}_p}{2I} = \mathbf{S}_q\) for some \(q\).
  Furthermore, since \(q_i = \frac{p_i}{2}\) and \(p_i \not = 0\),
  \(q_i \not = p_i\).
  Since \(\mathbf{S}_q\) and \(\mathbf{S}_p\) are equivalent,
  \(\mathbf{S}\) is not identifiable for \(p\).
\end{proof}

To complete the proof, we will need an additional lemma.
This lemma states that if a property \(P\), such as being not
identifiable, holds generically, then its complement, such as being
identifiable, does not hold generically.
The proof of this result is deferred to an appendix since it is uses
some technical concepts.

\begin{lemma}
  \label{lem:generic-complement}
  Let \(P\) be generic for \(\mathbf{S}\).
  Then, \(\neg P\) is not generic for \(\mathbf{S}\).
\end{lemma}

The above lemmas are sufficient to prove our result.
Hence, we state the proof of Theorem
\ref{thm:structured-identifiability} below.

\begin{proof}[Proof of Theorem \ref{thm:structured-identifiability}]
  By Lemma \ref{lem:generic-non-identifiable}, \(\mathbf{S}\) is generically
  not identifiable.
  The result then follows from Lemma \ref{lem:generic-complement}.
\end{proof}

To prove Theorem \ref{thm:generic-minimal-is-existential}, we will use
the following Lemma, which states that the parameter vectors for which
a structured linear system is not minimal are always elements of some
variety.
Since this lemma is nearly sufficient to prove our result, we will
state its proof in full.

\begin{lemma}
  \label{lem:generic-minimal-poly}
  Let \(\mathbf{S}\) be a structured linear system.
  Then, there exists a variety \(V_f\) such that for every parameter
  vector \(p\), \(p \in V_f\) if and only if \(\mathbf{S}_p\) is not minimal.
\end{lemma}
\begin{proof}
  Let \(x\) be a vector consisting of \(n\) indeterminates
  \(x_1,x_2,\cdots,x_n\).
  Then, \(\mathbf{S}_x\) is a system whose system matrices are polynomial
  matrices in the variables \(x_i\).
  Let \(O\) be the observability matrix of \(\mathbf{S}_x\) and \(C\) its
  controllability matrix.
  Clearly, \(O\) and \(C\) are polynomial matrices in the variables
  \(x_i\).
  Let \(f_O\) be the sum of the squares of all maximal-order minors of
  \(O\), and similarly for \(f_c\) using minors of \(C\).
  Thus, for any parameter vector \(p\), \(f_O(p) = 0\) if and only if
  \(\mathbf{S}_p\) is not observable, and \(f_C(p) = 0\) if and only
  if \(\mathbf{S}_p\) is not controllable.
  It follows that \(f_o(p)f_p(p)\) is zero if and only if \(\mathbf{S}_p\) is
  not minimal.
  Thus, the variety \(V_{f_of_p}\) is the variety we require.
\end{proof}

The proof of Theorem \ref{thm:generic-minimal-is-existential} using
the above lemma is stated below.

\begin{proof}[Proof of Theorem \ref{thm:generic-minimal-is-existential}]
  By Lemma \ref{lem:generic-minimal-poly}, there exists a variety
  \(V_f\) such that \(p \in V_f\) if and only if \(\mathbf{S}_p\) is not
  minimal.
  Thus, \(V_f\) is proper if and only if there exists a parameter
  vector \(p\) such that \(\mathbf{S}_p\) is minimal.
  Therefore, if there exists a parameter vector \(p\) such that
  \(\mathbf{S}_p\) is minimal, \(V_f\) is a proper variety containing the
  vectors \(q\) for which \(\mathbf{S}_q\) is not minimal.
  Thus, \(S\) is then generically minimal.
  If \(S\) is generically minimal, \(V_f\) must be proper, and
  so there exists a vector \(p \not \in V_f\) such that
  \(\mathbf{S}_p\) is minimal.
\end{proof}

Below, we state the proof of Theorem \ref{thm:generic-minimality}.

\begin{proof}[Proof of Theorem \ref{thm:generic-minimality}]
  To prove necessity, assume \(\mathbf{S}\) is generically minimal.
  Then, there exists a proper variety \(V_f\) such that if \(\mathbf{S}_p\) is
  not minimal, the parameter vector \(p \in V_f\).
  Let \(q\) be a parameter vector such that \(\mathbf{S}_q\) is not
  controllable.
  Then, \(\mathbf{S}_q\) is not minimal, and so \(q \in V_f\).
  Thus, \(\mathbf{S}\) is generically controllable.
  By the same argument, \(\mathbf{S}\) is generically observable.

  To prove sufficiency, assume \(\mathbf{S}\) is generically controllable and
  generically observable.
  Then, there exist proper varieties \(V_c\) and \(V_o\) such that \(p
  \in V_c\) if \(\mathbf{S}_p\) is not controllable and \(p \in V_o\) if
  \(\mathbf{S}_p\) is not observable.
  Let \(V_m\) be the variety \(\{p \in \mathbb{R}^n | c(p)o(p) =
  0\}\).
  Then, \(p \in V_m\) if and only if \(p \in V_c\) or \(p \in V_o\).
  Thus, if \(\mathbf{S}_p\) is not minimal, either \(p \in V_c\) or \(p \in
  V_o\), and so \(p \in V_m\).
  Since the set of all polynomial functions over the real numbers in
  \(n\) indeterminates forms an integral domain, \(c \cdot o \not = 0\), and
  so the variety \(V_m\) is proper.
  Therefore, \(\mathbf{S}\) is generically minimal.
\end{proof}

To prove Theorem \ref{thm:generic-minimality-graph}, we first state
graph-theoretical conditions for generic controllability.
The following theorem states the conditions for generic
controllability given by Dion et al.\cite{dion-et-al}.
Dion et al.\cite{dion-et-al} note that similar results hold for
generic observability.
Unfortunately, we have not been able to find a suitable reference for
these results.

\begin{theorem}
  \label{thm:generic-controllable}
  A structured linear system \stdnamedsysstruct is generically
  controllable if and only if the following conditions hold:
  \begin{enumerate}
  \item Every state variable \(x_i\) is the end vertex of some
    \(U\)-rooted simple path in \(G(\mathbf{S})\).
  \item There exists a disjoint union of a \(U\)-rooted simple path family
    and a cycle family that covers all state vertices.
  \end{enumerate}
\end{theorem}

To derive conditions for generic observability, we will first show
that a structured linear system is generically observable if and only
if its dual system is generically controllable.
The dual system of a linear system and a structured linear system is
stated below.
To define the dual of a structured linear system, we also define the
transpose of a structured matrix.

\begin{definition}
  Let \stdnamedsys be a linear system.
  The dual system of \(S\) is the system
  \namedsys{\dual{S}}{A^T}{C^T}{B^T}{D^T}.
\end{definition}

\begin{definition}
  \label{def:structured-transpose}
  Let \(A\) be a structured matrix determined by \(T(A)\).
  The transpose \(A^T\) of \(A\) is the structured matrix determined
  by \(T(A^T) = \{ (j,i) | (i,j) \in T(A) \}\).
\end{definition}

\begin{definition}
  Let \stdnamedsysstruct be a structured linear system.
  The dual system of \(\mathbf{S}\) is the structured linear system
  \namedsysstruct{\dualstruct{\mathbf{S}}}{A^T}{C^T}{B^T}{D^T}.
  Here, the transpose of a structured matrix defined in Definition
  \ref{def:structured-transpose} is used.
\end{definition}

\begin{remark}
  Let \(\mathbf{S}\) be a structured linear system.
  Then, for every system \(\mathbf{S}_p \in \mathbf{S}\), there exists a system \(\dualstruct{\mathbf{S}}_q
  \in \dualstruct{\mathbf{S}}\) that is the dual system of \(\mathbf{S}_p\).
  However, due to the way we have defined the standard
  parameterization of a structured linear system, the vector \(q\) is
  a permutation of the vector \(p\).
\end{remark}

The result that a structured linear system is generically observable
if and only if its dual is generically controllable is stated below as
a lemma.
The proof of this result uses the well-known fact that a linear system
is observable if and only if its dual is controllable and is deferred
to an appendix.

\begin{lemma}
  \label{lem:generic-observable-dual-controllable}
  A structured linear system \(\mathbf{S}\) is generically observable
  if and only if its dual \dualstruct{\mathbf{S}} is generically controllable.
\end{lemma}

To use the above result to derive graph-theoretical conditions for
generic observability, we need to formally state the relation between
the graph of a structured linear system and the graph of its dual.
We state this relation in the lemma below.
Intuitively, this lemma states that the graph of a dual system is
obtained from that of the original system by interchanging the inputs
and outputs of the system and reversing each arc of the original
system's graph.
The proof of this lemma is straightforward but somewhat involved and
is hence deferred to the appendix.

\begin{lemma}
  \label{lem:primal-dual-mapping}
  Let \(\mathbf{S}\) be a structured linear system.
  Furthermore, let \(G(\dualstruct{\mathbf{S}}) = (V_D,E_D)\) and
  \(G(\mathbf{S}) = (V_S,E_S)\).
  Then, the function \(f : V_D \rightarrow V_S\) given below has the
  following properties:
  \begin{enumerate}
  \item \(f\) is a bijection
  \item For all \(v_1,v_2 \in V_D\), \((v_1,v_2) \in E_D\) if and only if \((f(v_2),f(v_1))\in E_S\)
  \end{enumerate}

  \[ f(v) = 
  \begin{cases}
    x_i(G(\mathbf{S})) & \text{ if } v = x_i(G(\dualstruct{\mathbf{S}})) \\
    y_i(G(\mathbf{S})) & \text{ if } v = u_i(G(\dualstruct{\mathbf{S}})) \\
    u_i(G(\mathbf{S})) & \text{ if } v = y_i(G(\dualstruct{\mathbf{S}}))
  \end{cases}\]
\end{lemma}

In the lemmas below, we state the conditions that a system's graph
will satisfy if and only if the graph of the system's dual satisfies
the conditions for generic controllability.
Using Lemma \ref{lem:generic-observable-dual-controllable}, it is
clear that these conditions are a graph-theoretical characterization of
generic observability.
The proofs of these lemmas involve an intuitive application of Lemma
\ref{lem:primal-dual-mapping}.
The details of these proofs are deferred to an appendix.

\begin{lemma}
  \label{lem:y-topped-paths}
  Every state variable \(x_i(G(\mathbf{S}))\) is the first vertex of a
  \(Y\)-topped path if and only if every state variable
  \(x_i(G(\dualstruct{\mathbf{S}}))\) is the end vertex of a \(U\)-rooted path.
\end{lemma}

\begin{lemma}
  \label{lem:disjoint-unions}
  There exists a disjoint union of a \(Y\)-topped path family and a
  cycle family in \(G(\mathbf{S})\) that covers every vertex \(x_i\) if and
  only if a disjoint union of a \(U\)-rooted path family and a cycle
  family exists in \(G(\dualstruct{\mathbf{S}})\) that covers every vertex \(x_i\).
\end{lemma}

Using the above Lemmas, we can state the graph-theoretical
characterization of generic observability.
The proof of this result is straightforward and is hence deferred to
an appendix.

\begin{lemma}
  \label{lem:generic-observable}
  \(\mathbf{S}\) is generically observable if and only if every vertex \(x_i\) in
  \(G(\mathbf{S})\) is the first vertex of a \(Y\)-topped path and there exists
  a disjoint union of a \(Y\)-topped path family and a cycle family in
  \(G(\mathbf{S})\) that covers every vertex \(x_i\).
\end{lemma}

Below, we complete the proof of Theorem
\ref{thm:generic-minimality-graph}.

\begin{proof}[Proof of Theorem \ref{thm:generic-minimality-graph}]
  By Theorem \ref{thm:generic-minimality}, \(\mathbf{S}\) is generically minimal
  if and only if it is generically controllable and generically
  observable.
  Using Theorem \ref{thm:generic-controllable}, we find the first two
  conditions.
  The remaining conditions follow from Lemma \ref{lem:generic-observable}.
\end{proof}

Finally, we will show that the graph-theoretical conditions we have
previously derived are a necessary condition for a given linear system
to be minimal.
We state the proof of this result below.
Before we state this proof, we formally define the structured linear
system corresponding to the graph \(G(S)\) of a linear system \(S\).

\begin{definition}
  Let \stdnamedsys be a linear system.
  Furthermore, let \(T(A) = \{(i,j) | A_{ij} = 0\}\) and similarly for
  \(T(B)\),\(T(C)\) and \(T(D)\).
  Then, let \(A_s\) be determined by \(T(A)\) and similarly for
  \(B_s\), \(C_s\) and \(D_s\).
  The structured linear system corresponding to the graph of \(S\) is
  the system \namedsysstruct{\mathbf{S_s}}{A_s}{B_s}{C_s}{D_s}.
\end{definition}

\begin{proof}[Proof of Theorem \ref{thm:minimality-graph-necessary}]
Let \(S\) be a minimal linear system and \(\mathbf{S}'\) the structured linear
system corresponding to the graph \(G(S)\).
Since every entry of \(S\) is zero if and only if the corresponding
entry in \(\mathbf{S}'\) is a fixed zero, it is clear that the graphs \(G(S)\)
and \(G(\mathbf{S}')\) are identical.
Furthermore, it is clear that \(S \in \mathbf{S}'\).
Thus, since \(S\) is minimal, \(\mathbf{S}'\) is generically minimal.
Then, the graph \(G(\mathbf{S}')\) must satisfy the conditions of Theorem
\ref{thm:generic-minimality-graph}.
Since the graphs \(G(S)\) and \(G(\mathbf{S}')\) are identical, \(G(S)\) must
also satisfy these conditions, as claimed.
\end{proof}

\section{Conclusions}
In our introduction, we described a number of studies in which
researchers attempted to identify the structure of a dynamical system
from input/output data.
In this paper, we have considered to what degree this structure is
determined by the input/output relation of a linear system.

We began by applying linear transformations to systems.
As we saw in Subsection \ref{subsec:main-graph-isomorphism}, such
linear transformations may change the system's graph structure.
Furthermore, as we saw in Subsection \ref{subsec:main-cg-isomorphism}, even the
condensed graph of a system is not necessarily conserved by linear
transformations.
Finally, as we stated in Subsection \ref{subsec:main-cg-order}, even
the number of completely disconnected components in the system's graph
is not determined by input/output behavior.

The results described above indicate that many aspects of a system's
graph structure are not determined by its input/output relation.
This implies that to identify the structure of a linear system from
input/output data, we require assumptions about the system.
We have considered two possible forms of such assumptions.

The first kind of assumption we considered was that the system under
consideration has a particular canonical form.
For instance, the system might be a minimal SISO system with a
diagonal \(A\)-matrix.
We showed in Subsection \ref{subsec:main-cg-isomorphism} that we can
characterize the existence of (CG)-isomorphisms between two systems
satisfying this kind of assumption.
In addition, we showed in Subsection
\ref{subsec:main-structured-generic} that there exist certain
conditions that must be satisfied by the graph of a minimal system.

The second kind of assumption we considered was that the system was a
member of a particular structured system.
That is, we assumed that some edges could not occur in the system's
graph.
As we showed in Subsection \ref{subsec:main-structured-generic},
assumptions of this kind are not sufficient to uniquely identify a
system's parameters using input/output data.

To summarize, our results have two implications for the identification
of a system's structure.
First, a system's structure is not uniquely determined by the system's
input/output relation.
Thus, we require additional assumptions to identify a system's
structure.
Second, the assumption that a given set of edges does not occur in a
system's graph may be insufficient to identify a system's structure.
At the very least, this assumption is insufficient to uniquely
identify the system parameters.
Thus, even if the structure can be identified uniquely, the strength
of the influence of one variable on another cannot be quantified.

\subsection{Future work}
As we remarked in the introduction, we conjecture that some of our
results may apply, possibly in a modified form, to the models
considered by Hollanders\cite{hollanders} and Friston et
al.\cite{friston}.
We also remarked that the relation between the graph
structure of a vector AR model and the Granger causality criteria of
Goebel et al.\cite{goebel} is still unclear.
Thus, more work is needed to examine the implications of this paper
for the models considered by Hollanders, Friston et al.\ and Goebel et al.
\appendix
\section{Appendix}
This appendix consists of background material and technical proofs.
In the first subsection, we briefly recall the concepts of the
invariant polynomials and elementary divisors of a matrix.
The remaining subsections each contain the technical proofs that were
omitted from the corresponding section in the main text.
In addition, these subsections contain minor lemmas that are only
required for the technical proofs.

\subsection{Elementary divisors and invariant polynomials}
Consider an \(n \times n\) real matrix \(A\) and its characteristic matrix \charmat{A}.
Gantmacher \cite[Ch.~6]{gantmacher} shows that this characteristic matrix can
be transformed to a diagonal matrix D, as shown below, by elementary
row and column operations.
The diagonal elements of \(D\) are called the invariant polynomials of
the characteristic matrix \charmat{A} or equivalently those of the
matrix \(A\).
An important property of the polynomials \( i_1, i_2, i_3, \cdots
,i_n\) is that each divides the preceding one, that is \(i_i =
i_{i+1}p\), for some polynomial \(p\).

\[ D = \left[\begin{array}{cccc}
i_n & 0 & \cdots & 0 \\
0 & i_{n-1} & \cdots & 0 \\
\cdots & \cdots & \cdots & \cdots \\
0 & 0 & \cdots & i_1 
\end{array}
\right]\]

Gantmacher also shows an equivalent definition of these polynomials
using minors of \charmat{A}. Let \(D_i\) denote the greatest
common divisor of the minors of order \(i\) of \charmat{A}, with \(D_0
= 1\).
Then, we can equivalently define the invariant polynomials as ratios
of these greatest common divisors, as follows:
\[ i_1 = \frac{D_n}{D_{n-1}}, i_2 = \frac{D_{n-1}}{D_{n-2}}, \cdots,
i_n = \frac{D_1}{D_0}\]

In the above definition of \(i_1\), notice that \(D_n =
\dcharmat{A}\). This implies that \(i_1\) is the minimal polynomial of
A. As defined by Gantmacher \cite[Ch.~4]{gantmacher}, this polynomial
is the polynomial \(\psi\) of least degree such that \(\psi(A) = 0\),
where the first coefficient of \(\psi\) is taken to be 1.

Using the above definitions of the invariant polynomials, we can
define the elementary divisors of \(A\).
To do this, we factor each of the invariant polynomials \(i_j\) of
\(A\) into powers of irreducible polynomials \(\phi_i\), as shown below.
The powers of these polynomials \(\phi_i\) with exponents not equal to
zero that appear in this factorization are called the elementary
divisors of \(A\).

\begin{align*}
i_1 &= (\phi_1)^{k_{11}}(\phi_2)^{k_{12}}\cdots(\phi_m)^{k_{1m}}\\
i_2 &= (\phi_1)^{k_{21}}(\phi_2)^{k_{22}}\cdots(\phi_m)^{k_{2m}}\\
\vdots &\\
i_n &= (\phi_1)^{k_{n1}}(\phi_2)^{k_{n2}}\cdots(\phi_m)^{k_{nm}}\\
\end{align*}

\subsection{Technical proofs: Systems, graph structures and equivalent
  structures}
\label{app:prelims}
\begin{proof}[Proof of Observation \ref{obs:input-output-singletons}]
  For input vertices, it is clear that no path from a vertex \(v\) to
  the input \(u_i(G(S))\) exists, as input vertices have indegree
  zero.
  Thus, there exist no vertices \(v\) other than \(u_i(G(S))\) such
  that \(v \leftrightarrow u_i(G(S))\).
  Therefore, \(u_i(G(S))\) is the sole member of its component.
  The proof for \(y_i(G(S))\) is similar.
\end{proof}
\subsection{Technical proofs: Graph isomorphism and its inadequacy}
\begin{proof}[Proof of Lemma \ref{lem:tr-isomorphism-equiv}]
To prove the first condition, consider the identity function on
\(V(G(S))\).
Clearly, the identity is bijective and satisfies the second condition
of Definition \ref{def:graph-isomorphism}.
It is also clear that the identity is type-restricted, and so \(S
\simeq S\).

To prove symmetry, let \(S \simeq S'\). Then, there exists a
type-restricted isomorphism \(\phi : V(G(S)) \rightarrow V(G(S'))\).
Clearly, the inverse \(\phi^{-1}\) is bijective.
This inverse also satisfies the second condition of Definition
\ref{def:graph-isomorphism}, since \(\phi\) is surjective.
For the same reason, the type-restriction condition must also be
satisfied.
Therefore, \(S' \simeq S\).

To complete the proof, let \(S \simeq S'\) and \(S' \simeq S''\).
Then, there exist type-restricted isomorphisms \(\phi_1 : V(G(S))
\rightarrow V(G(S'))\) and \(\phi_2 : V(G(S')) \rightarrow
V(G(S''))\).
The composition \(\phi = \phi_2 \circ \phi_1 \) is then a bijection
from \(V(G(S))\) to \(V(G(S''))\).
This composition also satisfies the other conditions for a
type-restricted isomorphism, as can be readily verified.
\end{proof}

\begin{proof}[Proof of Lemma \ref{lem:permutation-orthogonal}]
  Since the rows of \(P(e_1,e_2,\cdots,e_n)\) are a permutation of the rows of \(I_n\),
  they are clearly an orthogonal set of unit vectors.
  Therefore, \(P(e_1,e_2,\cdots,e_n)\) is clearly orthogonal.
\end{proof}

\subsection{Technical proofs: Condensed-graph isomorphism}
\begin{proof}[Proof of Lemma \ref{lem:trap-non-observability}]
Let \(v_1,v_2,\cdots,v_m\) be the indices of the state variables in
the trap in the graph of \(G(S)\).
Then, it is clear that \((x_{v_i},y_j)\) is not an edge of \(G(S)\) for
all \(i\) and \(j\), or equivalently, that \(C_{jv_i} = 0\).
For the same reason, for all integers \(j\) such that \(j \not = v_i\)
for all \(i\), \(A_{jv_i} = 0\) for all \(i\).
Consider \((CA)_{jv_i} = \sum_k C_{jk}A_{kv_i}\).
From the above discussion, we have for all \(k\) that either \(C_{jk}
= 0\) or \(A_{kv_i} = 0\).
Thus, \((CA)_{jv_i} = 0\).
By repeating this argument, we can show that the same holds for
\(CA^k\) for all \(k\).
Therefore, the columns \(v_1,v_2,\cdots,v_m\) of \(\mathfrak{O}\) are
zero.
Since \(\mathfrak{O}\) has \(n\) columns, \(\mathfrak{O}\) cannot be
of rank \(n\).
\end{proof}

\begin{proof}[Proof of Lemma \ref{lem:minpoly-non-minimality}]
Let \(i_1 \not = \dcharmat{A}\).
Since \(\dcharmat{A} = i_1 \gcd(\cadj{A})\), where \cadj{A} is the
adjugate of \charmat{A}, \(\gcd(\cadj{A}) \not = 1\).
Consider the nominal transfer function \( H =
\frac{C\cadj{A}B}{\dcharmat{A}}\).
Since \(\dcharmat{A} = \gcd(\cadj{A})i_1\) and \(\cadj{A} =
\gcd(\cadj{A})\Gamma\), for some matrix \(\Gamma\), \( H =
\frac{\gcd(\cadj{A})C\Gamma B}{\gcd(\cadj{A})i_1} = \frac{C\Gamma
  B}{i_1}\).
Thus, \(H \) has a pole-zero cancellation and so the realization
\((A,B,C,D)\) is non-minimal.
\end{proof}

\begin{proof}[Proof of Lemma \ref{lem:elem-div-structure}]
Since \(\phi\) is irreducible, either \(\phi(0) \not = 0\) or
\(\phi(\lambda) = \lambda^m\) for some \(m \geq 1\).
In the former case, this implies that in \(\phi(\lambda) = \lambda^m +
a_1\lambda^{m-1} + \cdots + a_m\), \(a_m \not = 0\).
Therefore, the top-right entry in the \(n \times n \) companion matrix \(L\)
corresponding to \(\phi(\lambda)^k\) is non-zero.
Then, an edge  \((x_n,x_1)\) exists in the graph of \(L\).
Since the fixed elements equal to 1 on \(L\)'s subdiagonal correspond
to edges \((x_i,x_{i+1})\) for \(1 \leq i < n\), this graph has a
Hamiltonian cycle.

Otherwise, if \(\phi(\lambda) = \lambda^m\), \(\phi(\lambda)^k =
\lambda^{mk}\).
Thus, in the \(n \times n \) companion matrix \(L'\) corresponding to this case, the
last column is zero.
Then, the only edges that exist in the graph of \(L'\) are the edges
\((x_i,x_{i+1})\) for \(1 \leq i < n\), implying that this graph is a
directed path.
\end{proof}

\begin{proof}[Proof of Lemma \ref{lem:second-natural-component-conditions}]
Since the realization \(S\) is minimal, every vertex in every
component cannot be part of either an unreachable set or a trap.
In the first case, if the component is Hamiltonian, if no edge of the
form \((u,x_i)\) exists for \(x_i\) in the component, the entire
component is unreachable.
Similarly, if no edge of the form \((x_j,y)\) exists for \(x_j\) in
the component, the entire component is a trap.
Therefore, edges of these forms must exist.
Since the component is Hamiltonian, these conditions are also
sufficient for the component to contain neither traps nor unreachable
sets.

In the second case, the component is a directed path consisting of the
vertices \(x_i,x_{i+1},\cdots,x_{i+m}\).
If the edge \((u,x_i)\) does not exist in \(G(S)\), \(\{x_i\}\) is an
unreachable set, and so \(S\) is non-minimal.
Similarly, if \((x_{i+m},y)\) does not exist in \(G(S)\),
\(\{x_{i+m}\}\) is a trap, and so \(S\) is non-minimal.
Thus, the stated conditions must hold.
Furthermore, if the conditions are satisfied, paths from \(u\) to
\(x_i\) and from \(x_i\) to \(y\) exist in \(G(S)\) for all \(x_i\) in
the component, and so the component does not contain traps or
unreachable sets.

By Lemma \ref{lem:elem-div-structure}, it is clear that each component
consisting of state variables in \(G(S)\) is covered by one of the
two cases above.
Therefore, no such component can contain traps or unreachable sets,
and so \(G(S)\) can contain neither traps nor unreachable sets.
\end{proof}

\begin{proof}[Proof of Lemma \ref{lem:minimal-irreducible-poly-uniqueness}]
Let \(A\) have two elementary divisors \(\phi(\lambda)^{k_1}\) and
\(\phi(\lambda)^{k_2}\).
Since these elementary divisors are powers of the same irreducible
polynomial, they cannot occur in the same invariant polynomial.
After all, if they did, this invariant polynomial would only have a
single elementary divisor \(\phi(\lambda)^{k_1+k_2}\).
Therefore, two or more invariant polynomials of \(A\) are not
equal to 1.
Then, the minimal polynomial \(i_1\) of \(A\) cannot coincide with
\dcharmat{A}, and so by Lemma \ref{lem:minpoly-non-minimality}, \(S\)
is non-minimal.
\end{proof}

\subsection{Technical proofs: Components of condensed graphs}
\label{app:tech-condensed-components}
\begin{lemma}
\label{lem:tuple-realization}
Let \(\phi_i\), \(1 \leq i \leq l\) be the irreducible polynomials
that divide \dcharmat{A}.
Furthermore, let \(E_i\), \(1 \leq i \leq m\) be sets of elementary
divisors of \(A\) such that no two elements of \(E_i\) are of the form
\(\phi_j^l\) for the same \(j\), for all \(i\).
Additionally, let each elementary divisor of \(A\) be an element of
exactly one set \(E_i\).
Then, the \(m\)-tuple \((E_1,E_2,\cdots,E_m)\) corresponds to a
block-companion realization similar to \stdnamedsys.
\end{lemma} 
\begin{proof}
Let \(l_i = \prod_{e \in E_i}e\),\(1 \leq i \leq m\), and \(L_i\) be the companion matrix
corresponding to \(l_i\).
We claim that \(A' = \{L_i\}\) is similar to A.
Since \(A'\) is block-diagonal, the elementary divisors of \(A'\) are
those of the matrices \(L_i\), by Theorem \ref{thm:elem-div-diagonal}.
Since each of the matrices \(L_i\) is a companion matrix corresponding
to the product of the elements of \(E_i\), the elementary divisors of
\(L_i\) are the elements of \(E_i\).
Thus, each elementary divisor of \(A\) is an elementary divisor of
\(A'\), since we require that each such divisor is an element of some
\(E_i\).
In addition, \(A'\) can have no other elementary divisors, as each
element of each set \(E_i\) is an elementary divisor of \(A\) and each
elementary divisor of \(A\) is in exactly one set \(E_i\).
Therefore, \(A'\) and \(A\) are similar.
Thus, there exists a matrix \(T\) such that \(A' = TAT^{-1}\).
Therefore, \transys{S}{T} is a block-companion realization similar to
\(S\), as claimed.
\end{proof}

\begin{proof}[Proof of Lemma \ref{lem:diag-block-lb-sharpness}]
Let \(\phi_i\),\(1 \leq i \leq m\) be the irreducible polynomials that
divide \dcharmat{A}. 
Furthermore, let \(e_{ij}\) be the \(j\)-th elementary divisor of
\(A\) of the form \(\phi_i^l\), in some arbitrary order.
Let \(E_{ij}\) be a set of elementary divisors of \(A\), defined as
follows:
\[ E_{ij} = \begin{cases}
\{e_{ij}\} &\text{if } e_{ij} \text{ exists} \\
\emptyset & \text{otherwise}
\end{cases}\]
Clearly, each elementary divisor \(e_{ij}\) is an element of only
\(E_{ij}\) and no other set \(E_{lm}\).
Thus, the \(k\)-tuple \(t = (\bigcup_{1 \leq i \leq m}E_{i1},\bigcup_{1 \leq i
  \leq m}E_{i2},\cdots,\bigcup_{1 \leq i \leq m}E_{ik})\) consists
of sets of elementary divisors of \(A\).
It is clear that each elementary divisor of \(A\) is an element
of exactly one element of the tuple \(t\).
Furthermore, each element of the tuple \(t\) contains at most 1
elementary divisor of the form \(\phi_i^l\) for each polynomial
\(\phi_i\).
Thus, by Lemma \ref{lem:tuple-realization}, \(t\) corresponds to a
realization \(S'\) similar to \(S\).
\(S'\) is a block-companion realization with \(k\) diagonal blocks, as claimed.
\end{proof}

\begin{proof}[Proof of Lemma \ref{lem:diag-block-ub-sharpness}]
Let \(A'\) be the second natural normal form of \(A\).
Then, \(A'\) has \(n\) diagonal blocks and is a block-companion
matrix.
Furthermore, \(A'\) is similar to \(A\).
Thus, a matrix \(T\) exists such that \(A' = TAT^{-1}\).
Therefore, \transys{S}{T} is the required realization.
\end{proof}

\begin{proof}[Proof of Lemma \ref{lem:diag-block-inc}]
Let \(t = (E_i)\) be an \(l\)-tuple of sets of elementary divisors,
where \(E_i\) consists of the elementary divisors of the \(i\)-th
diagonal block of \(A\).
Assume \(l \not = n\).
Then, since each elementary divisor of \(A\) is a member of some set
\(E_i\), at least one set \(E_i\) consists of two or more elements.
Let \(j\) be an integer such that \(E_j\) consists of two or more
elements and select an arbitrary element \(e\) of \(E_j\).
We claim that \(t' = (E_1,\cdots,E_{j-1},E_j
\backslash \{e\},E_{j+1},\cdots,E_l,\{e\})\) is an \((l+1)\)-tuple satisfying
the conditions of Lemma \ref{lem:tuple-realization}.
Clearly, since none of the sets \(E_i\) contains two elementary
divisors of the form \(\phi_i^l\) for some irreducible polynomial
\(\phi_i\), neither do the elements of \(t'\).
The elements of \(t'\) are also clearly sets of elementary divisors of
\(A\).
Furthermore, since each elementary divisor is an element of exactly
one set \(E_i\), the same holds for the elements of \(t'\).
Thus, by Lemma \ref{lem:tuple-realization}, the tuple \(t'\)
corresponds to a block-companion realization similar to \(S\) with
\(l+1\) diagonal blocks.
\end{proof}
\subsection{Technical proofs: Structured systems and their graphs}
\begin{proof}[Proof of Lemma \ref{lem:trans-equivalence}]
  Note that \(\transys{S}{M} = (MAM^{-1},MB,CM^{-1},D)\).
  First, we show by induction that \(x_{\transys{S}{M},k} = Mx_{S,k}\).
  Since we use zero initial conditions, \(x_{\transys{S}{M},k} = 0 =
  M0 = Mx_{S,0}\).
  Inductively, \(x_{\transys{S}{M},k+1} =
  MAM^{-1}x_{\transys{S}{M},k}+MBu_k = M(AM^{-1}Mx_{S,k}+Bu_k) =
  M(Ax_{S,k} + Bu_k) = Mx_{S,k+1}\).
  Therefore, \(y_{\transys{S}{M},k} = CM^{-1}x_{\transys{S}{M},k}+Du_k
  = Cx_{S,k}+Du_k = y_{S,k}\).
\end{proof}

\begin{proof}[Proof of Lemma \ref{lem:generic-complement}]
  Since \(P\) is generic for \(\mathbf{S}\), there exists a proper variety \(V_f\)
  such that \(\{p \in \mathbb{R}^n | P \text{ does not hold for } \mathbf{S}_p\}
  \subset V_f\).
  Suppose \(\neg P\) is generic for \(\mathbf{S}\), i.e. there exists a proper variety
  \(V_g\) such that \(\{p \in \mathbb{R}^n | P \text{ holds for } \mathbf{S}_p\} \subset V_g\).
  Then, \(\mathbb{R}^n \subset V_f \cup V_g\) and so \(V_f \cup V_g
  = \mathbb{R}^n\).
  Let \(h(x) = f(x)g(x)\).
  Then, \(V_f \cup V_g = V_h = \mathbb{R}^n\).
  Therefore, \(h\) is the zero function.
  But then, since the domain of polynomials over \(\mathbb{R}\) in
  \(n\) indeterminates is an integral domain, either \(f\) or \(g\)
  must be the zero function.
  Since \(V_f\) is proper, \(f\) is non-zero for some \(x\), and so
  \(g\) must be the zero function.
  But then, either \(V_g\) is not proper, which is impossible.
\end{proof}

\begin{lemma}
  Let \stdnamedsys be a linear system.
  Then, \(S\) is observable if and only if \dual{S} is controllable.
\end{lemma}
\begin{proof}
  This result is implicitly given in the textbook by Kailath
  \cite{kailath}.
  To formally prove it, notice that
  \(\namedsys{\dual{S}}{A^T}{C^T}{B^T}{D^T}\).
  Thus, the controllability matrix of \(\dual{S}\) is \(C_D =
  \begin{bmatrix}
    C^T &A^TC^T &\cdots &(A^T)^{n-1}C^T
  \end{bmatrix}\).
  Since the observability matrix of \(S\) is given by \(O_S = 
  \begin{bmatrix}
    C \\
    CA\\
    \vdots \\
    CA^{n-1}
  \end{bmatrix}\), we notice that \(C_D = O_S^T\).
  Thus, these matrices have the same rank, completing our proof.
\end{proof}

\begin{proof}[Proof of Lemma \ref{lem:generic-observable-dual-controllable}]
  To prove necessity, assume \(\mathbf{S}\) is generically observable.
  Thus, there exists a proper variety \(V_f\) such that if \(\mathbf{S}_p\) is
  not observable, \(p \in V_f\).
  Let \(g\) be the polynomial obtained from \(f\) by permuting the
  indeterminates in \(f\) such that for all \(p\) and \(q\) such that \(\dual{\mathbf{S}_p} =
  \dualstruct{\mathbf{S}}_q\), \(f(p) = 0\) if and only if \(g(q) = 0\).
  Let \(q\) be an arbitrary parameter vector such that
  \(\dualstruct{\mathbf{S}}_q\) is not controllable and let \(p\) be the vector
  such that \(\dual{\mathbf{S}_p} = \dualstruct{\mathbf{S}}_q\).
  Then, \(\mathbf{S}_p\) is not observable, and so \(p \in V_f\).
  But then, \(q \in V_g\).
  Since \(f\) is not identically zero, neither is \(g\), and so
  \(V_g\) is proper.
  Thus, \(\dualstruct{\mathbf{S}}\) is generically controllable.
  
  The sufficiency of the condition follows from a similar argument.
\end{proof}

\begin{proof}[Proof of Lemma \ref{lem:primal-dual-mapping}]
  The first property, that \(f\) is a bijection, is clear from the
  definition of \(f\) and the definition of the dual \dualstruct{\mathbf{S}}.
  
  To prove the second property, let \stdnamedsysstruct.
  Then, \namedsysstruct{\dualstruct{\mathbf{S}}}{A^T}{C^T}{B^T}{D^T}.
  We will consider all the types of
  edges that occur in \(G(\dualstruct{\mathbf{S}})\).
  First, consider edges of the form \((x_i,x_j)\), which exist in
  \(G(\dualstruct{\mathbf{S}})\) if and only if \(A^T_{ji}\) is not a fixed
  zero.
  Since \(A^T_{ji}\) is a fixed zero if and only if \(A_{ij}\) is a
  fixed zero, this edge exists if and only if \((x_j,x_i)\) is an edge
  of \(G(\mathbf{S})\).
  Since \(f(x_i) = x_i\), this proves our condition for edges of this
  form.
  
  Second, consider edges of the form \((u_i,x_j)\), which exist in
  \(G(\dualstruct{\mathbf{S}})\) if and only if \(C^T_{ji}\) is not a fixed
  zero.
  Since \(C^T_{ji}\) is a fixed zero if and only if \(C_{ij}\) is a
  fixed zero, this edge exists if and only if \((x_j,y_i)\) is an edge
  of \(G(\mathbf{S})\).
  Since \(f(u_i) = y_i\) and \(f(x_j) = x_j\), this proves our
  condition for edges of this type.

  The arguments for the remaining edges are similar.
\end{proof}

\begin{proof}[Proof of Lemma \ref{lem:y-topped-paths}]
  To prove necessity, let \(x_i v_1 v_2 \cdots v_k y_j\) be a
  \(Y\)-topped path in \(G(\mathbf{S})\).
  Let \(f\) be the mapping of Lemma \ref{lem:primal-dual-mapping}.
  Then, since \(f\) is a bijection and the vertices \(x_i\), \(v_i\)
  and \(y_j\) are all distinct, so are the vertices
  \(f^{-1}(x_i)\),\(f^{-1}(v_i)\) and \(f^{-1}(y_j)\).
  Furthermore, since edges \((x_i,v_1)\), \((v_i,v_{i+1})\) and
  \((v_k,y_j)\) exist in \(G(\mathbf{S})\), edges
  \((f^{-1}(v_1),f^{-1}(x_i))\), \((f^{-1}(v_{i+1}),f^{-1}(v_i))\) and
  \((f^{-1}(y_j),f^{-1}(v_k))\) exist in \(G(\dualstruct{\mathbf{S}})\).
  Thus, using the definition of \(f\), we find that \(u_j v'_1 v'_2
  \cdots v'_k x_i\) is a path in \(G(\dualstruct{\mathbf{S}})\), where \(v'_i = f^{-1}(v_i)\).
  
  The sufficiency of the condition follows from the same argument.
\end{proof}

\begin{proof}[Proof of Lemma \ref{lem:disjoint-unions}]
  To prove necessity, suppose a disjoint union as described above exists.
  Then, every state variable \(x_i\) is covered by either a
  \(Y\)-topped path or a cycle.
  Using the proof of Lemma \ref{lem:y-topped-paths}, a state variable
  covered by a \(Y\)-topped path is covered by a \(U\)-rooted path in
  \(G(\dualstruct{\mathbf{S}})\).
  A similar argument shows that a state variable covered by a cycle in
  \(G(\mathbf{S})\) will also be covered by a cycle in \(G(\dualstruct{\mathbf{S}})\).
  Thus, every \(Y\)-topped path in the path family corresponds to a
  \(U\)-rooted path in \(G(\dualstruct{\mathbf{S}})\), and every cycle
  corresponds to a cycle.
  It remains to be shown that these paths and cycles are mutually
  disjoint.
  Thus, let \(v_1v_2\cdots v_k\) and \(v'_1v'_2\cdots v'_k\) be
  disjoint paths or cycles in \(G(\mathbf{S})\).
  Then, the corresponding paths or cycles in \(G(\dualstruct{\mathbf{S}})\) are
  given by \(f^{-1}(v_k)\cdots f^{-1}(v_2)f^{-1}(v_1)\) and
  \(f^{-1}(v'_k)\cdots f^{-1}(v'_2)f^{-1}(v'_1)\).
  Suppose these paths are not disjoint.
  Then, there exist \(i\) and \(j\) such that \(f^{-1}(v_i) =
  f^{-1}(v'_j)\).
  Since \(f^{-1}\) is injective, this implies that \(v_i = v'_j\).
  This contradicts our assumption that the paths \(v_1v_2 \cdots v_k\)
  and \(v'_1 v'_2 \cdots v'_k\) were disjoint.
  Thus, the corresponding paths in \(G(\dualstruct{\mathbf{S}})\) must be
  disjoint.

  A similar argument in the other direction proves the sufficiency of
  the condition.
\end{proof}

\begin{proof}[Proof of Lemma \ref{lem:generic-observable}]
  By Lemma \ref{lem:generic-observable-dual-controllable}, \(\mathbf{S}\) is
  generically observable if and only if the dual \dualstruct{\mathbf{S}} is
  generically controllable.
  Furthermore, by Theorem \ref{thm:generic-controllable},
  \dualstruct{\mathbf{S}} is generically controllable if and only if every
  vertex \(x_i\) in \(G(\dualstruct{\mathbf{S}})\) is the end vertex of a
  \(U\)-rooted path and there exists a disjoint union of a
  \(U\)-rooted path family and a cycle family in \(G(\dualstruct{\mathbf{S}})\)
  that covers every vertex \(x_i\).
  By Lemmas \ref{lem:y-topped-paths} and \ref{lem:disjoint-unions},
  these conditions are equivalent to the stated conditions on
  \(G(\mathbf{S})\), completing the proof.
\end{proof}
\bibliography{refs}
\end{document}